\documentclass[11pt, notitlepage]{article}   

\usepackage{amsmath,amsthm,amsfonts}   

\usepackage{bbm} 
\usepackage{amscd, array}
\usepackage{amsfonts}
\usepackage{amssymb}
\usepackage{color}      
\usepackage{epsfig}
\usepackage{graphicx}           
\usepackage{algorithm}
\usepackage[noend]{algpseudocode}

\usepackage{tikz}
\usetikzlibrary{shapes}
\usetikzlibrary{matrix}
\usetikzlibrary{positioning}
\usetikzlibrary{fit}

\tikzstyle{res}=[circle,thick,minimum size=4mm,draw=black,fill=red,inner sep=1pt]
\tikzstyle{non-res}=[circle,thick,minimum size=4mm,draw=black,inner sep=1pt]
\tikzstyle{light-res}=[circle,thick,minimum size=4mm,draw=black,fill=red!40,inner sep=1pt]
\tikzstyle{blue}=[circle,thick,minimum size=4mm,draw=black,fill=blue!20,inner sep=1pt]

\newtheorem{theorem}{Theorem}[section]

\theoremstyle{remark}
\newtheorem*{remark}{Remark}
\theoremstyle{plain}

\newtheorem{thm}[theorem]{Theorem}

\newtheorem{problem}[theorem]{Problem}

\def\finf{\mathop{{\rm I}\kern -.27 em {\rm F}}\nolimits}

\usepackage{graphicx} 
\usepackage[margin=1.5in]{geometry}
\usepackage{amsmath,amsthm,amssymb}
\usepackage[hidelinks]{hyperref}
\usepackage{xcolor}
\usepackage{booktabs}
\usepackage{tikz}
\usepackage{enumitem}

\usepackage{mathtools}
\usepackage{url}
\newcommand{\sat}{\operatorname{sat}}
\newcommand{\ssat}{\operatorname{ssat}}
\newcommand{\ex}{\operatorname{ex}}

\newtheorem{lem}[theorem]{Lemma}
\newtheorem{cor}[theorem]{Corollary}
\newtheorem{conj}[theorem]{Conjecture}

\newtheorem{prop}[theorem]{Proposition}
\newtheorem{claim}[theorem]{Claim}
\theoremstyle{definition}
\newtheorem{definition}[theorem]{Definition}
\newtheorem{obs}[theorem]{Observation}
\newtheorem{example}[theorem]{Example}
\theoremstyle{remark}

\begin{document}

\title{Saturation of 0-1 Matrices}
\author{Andrew Brahms, Alan Duan, Jesse Geneson, Jacob Greene}
\date{February 28, 2025}
\maketitle

\begin{abstract}
    A 0-1 matrix $M$ \textit{contains} a 0-1 matrix $P$ if $M$ has a submatrix $P'$ which can be turned into $P$ by changing some of the ones to zeroes. Matrix $M$ is $P$-\textit{saturated} if $M$ does not contain $P$, but any matrix $M'$ derived from $M$ by changing a zero to a one must contain $P$. The \textit{saturation function} $\sat(n,P)$ is defined as the minimum number of ones of an $n \times n$ $P$-saturated 0-1 matrix. Fulek and Keszegh showed that each pattern $P$ has $\sat(n,P) = O(1)$ or $\sat(n,P) = \Theta(n)$. This leads to the natural problem of classifying forbidden 0-1 matrices according to whether they have linear or bounded saturation functions. Some progress has been made on this problem: multiple infinite families of matrices with bounded saturation function and other families with linear saturation function have been identified. We answer this question for all patterns with at most four ones, as well as several specific patterns with more ones, including multiple new infinite families. We also consider the effects of certain matrix operations, including the Kronecker product and insertion of empty rows and columns. Additionally, we consider the simpler case of fixing one dimension, extending results of (Fulek and Keszegh, 2021) and (Berendsohn, 2021). We also generalize some results to $d$-dimensional saturation.

\end{abstract}

\section{Introduction}

All matrices in this paper are 0-1 matrices. The \textit{weight} of a matrix $P$ is the number of 1-entries in $P$. A 0-1 matrix $M$ \textit{contains} $P$ if some submatrix $P'$ of $M$ can be turned into $P$ by changing some number of ones to zeroes. Otherwise we say that $M$ \textit{avoids} $P$.

A well-studied combinatorial problem is the \textit{extremal function} $\ex(m,n,P)$, which is the maximum weight of an $m\times n$ matrix $M$ which avoids $P$; particular interest has centered around the square case $\ex(n,P) = \ex(n,n,P)$. The impetus for this study originated with analyzing the time complexity of an algorithm for finding the shortest rectilinear path between two points in a grid with obstacles \cite{mitchell}. This analysis was done by Bienstock and Gy\"ori \cite{bienstock}, who showed that the complexity was bounded by $\ex(n,P)$ for certain sparse forbidden patterns $P$. 

The most well-known use of the extremal function so far has been for proving the Stanley-Wilf conjecture, which gave a bound on the number of permutations of length $n$ avoiding a given forbidden permutation $\pi$. Klazar \cite{klazar} reduced this to a problem about $\ex(n,P)$, which was resolved by Marcus and Tardos \cite{MT}.

Recently, interest has also grown around a new function for pattern avoidance in 0-1 matrices, which is our main object of study. We say that a matrix $M$ is \textit{$P$-saturated} if $M$ does \textit{not} contain $P$, but changing any zero in $M$ to a one introduces a copy of $P$; that is to say, for all matrices $M'$ derived from $M$ by changing any zero in $M$ to a one, $M'$ contains $P$.
\begin{definition}[Saturation Function]
    The \textit{saturation function $\sat(m,n,P)$} is the minimum weight of an $m \times n$ $P$-saturated 0-1 matrix; we write $\sat(n,n,P) = \sat(n,P)$.
\end{definition}

Note that the extremal function $\ex(n,P)$, on the other hand, is the \textit{maximum} weight of an $n \times n$ $P$-saturated 0-1 matrix.

The saturation function, along with other related questions, was first studied by Brualdi and Cao in \cite{BC}, although it was Fulek and Keszegh \cite{FK} who later introduced the terminology of $\sat(n,P)$. Brualdi and Cao found certain classes of matrices for which $\sat(n,P)$ and $\ex(n,P)$ coincide. This work was continued by Fulek and Keszegh, whose paper answers some of the questions asked by Brualdi and Cao. However, though these early results suggested that $\sat(n,P)$ and $\ex(n,P)$ might exhibit similar behavior, Fulek and Keszegh found a marked difference.

Fulek and Keszegh found that the saturation function exhibits a certain dichotomy in its asymptotic behavior: all patterns have saturation function either bounded or (asymptotically) linear in $n$. Fulek and Keszegh identified multiple infinite classes of matrices with linear saturation function, and later work by Geneson \cite{geneson21} and Berendsohn \cite{berendsohn21, berendsohn23} has identified numerous other such classes, as well as other infinite classes of matrices with bounded saturation function. In particular, the asymptotic behavior of the saturation function of permutation matrices was completely characterized in \cite{berendsohn23}. However, the full classification for all 0-1 matrices remains incomplete. Our paper contributes to this classification, as well as resolving some other related questions about saturation.

In Section~\ref{sec:pastresults}, we briefly describe some results about the extremal function, before discussing in depth the past results that we will need on saturation. In Section~\ref{sec:1dsat}, we present several new results to develop the theory of $\sat(m_0,n,P)$ and $\sat(m,n_0,P)$ for fixed $m_0$, $n_0$. In particular, Fulek and Keszegh asked whether there exists a finite algorithm to determine whether $\sat(n,P)$ is bounded; we answer an analogue of this question for $\sat(m_0,n,P)$ in the affirmative. Also, we fully characterize semisaturation with one dimension fixed, we offer a corrected proof for a result of Fulek and Keszegh relating $\sat(m_0,n,P)$ and $\sat(m_1,n,P)$ when $m_1>m_0$, and we generalize a result of Berendsohn relating $\sat(m_0,n,P)$ and $\sat(m,n_0,P)$ to $\sat(n,P)$ to a wider class of patterns $P$.

In Sections~\ref{sec:insert} and \ref{sec:kron}, we examine the relationship between saturation and various operations, including the Kronecker product and the effect of inserting empty rows and columns into a pattern. In Section~\ref{sec:ddim}, we prove some new results about saturation functions of multidimensional 0-1 matrices. In particular, we provide the first example of a $d$-dimensional 0-1 matrix with bounded saturation function for $d > 2$; we generalize a result of Geneson and Tsai \cite{gt23} to the case of fixing certain dimensions while the others go to infinity; and we derive a sharp upper bound on the value of the $d$-dimensional saturation function over all forbidden $d$-dimensional 0-1 matrices of dimensions $k_1 \times k_2 \times \dots \times k_d$.

In Section~\ref{sec:witness_graphs}, we present a new construction of a graph associated to a witness for a pattern. We find some general results about these graphs for all patterns and then apply them to classify the saturation of two specific patterns. In Section~\ref{sec:lin_prog}, we present a linear program that characterizes saturation, and we apply it to evaluate small values of the saturation function for some patterns. Finally, in Section~\ref{sec:patterns}, we classify the saturation function for all patterns with at most $4$ ones or dimensions at most $4 \times 4$, along with several infinite families of patterns and various other miscellaneous examples. We conclude in Section~\ref{sec:conc} by discussing some open problems and future research directions.

\section{Past Results}
\label{sec:pastresults}

We will represent 0-1 matrices using a $\bullet$ to represent a one, an empty space or a $\cdot$ to represent a zero, and a $\ast$ to represent either. In line with this notation, we call a row or column \textit{empty} if it contains no ones.

\subsection{The Extremal Function}

There has been much work in determining the growth rate of the extremal function for many classes of matrices. Note that the extremal function is always $O(n^2)$ since there are $n^2$ total entries in an $n\times n$ matrix. On the other hand, $\ex(n,P)$ is always $\Omega(n)$ if $P$ is not $1 \times 1$ or an all-zeroes matrix. If $Q$ contains $P$, then clearly we have $\ex(n, P) \le \ex(n, Q)$. Rotating a forbidden pattern $P$ ninety degrees or reflecting it vertically or horizontally does not change the value of $\ex(n, P)$. A number of operations have been found which only change the extremal function by at most a constant factor (see, e.g., \cite{tardos, keszegh, pettie11}). 

There are many examples of matrices $P$ with $\ex(n,P) = \Theta(n)$. For example, any all-ones matrix $P$ with a single row or a single column of length $k$ has $\ex(n,P)= (k-1)n$. As mentioned in the introduction, Marcus and Tardos \cite{MT} proved that every permutation matrix $P$ has $\ex(n, P) = O(n)$, thereby proving the Stanley-Wilf conjecture. Geneson \cite{geneson2009} strengthened the result by showing that every 0-1 matrix obtained by doubling the columns of a permutation matrix has linear extremal function. Fulek found a new family of 0-1 matrices with linear extremal function \cite{fulek} using a reduction to bar visibility graphs. Another family of 0-1 matrices with linear extremal functions was found in \cite{GS15} using a reduction to bar visibility hypergraphs. 

Growth rates strictly between $\Theta(n)$ and $\Theta(n^2)$ are also possible, including patterns with $\Theta(n\log n)$ which were used in the proof of the $O(n \log n)$ upper bound on the maximum number of unit distances in a convex $n$-gon \cite{ngon}. Keszegh used a connection between forbidden 0-1 matrices and forbidden sequences to show that every 0-1 matrix $P$ with a single one in every column has $\ex(n, P) = O(n2^{\alpha(n)^t})$ for some $t$, where $\alpha(n)$ is the incredibly slow-growing inverse Ackermann function \cite{keszegh}. Among the 0-1 matrices with at most four ones, the greatest order of magnitude of the extremal function is $\Theta(n^{3/2})$, attained by a $2 \times 2$ all-ones matrix. More generally, the extremal function of a $j \times k$ all-ones matrix is $O(n^{2-1/k})$ using an argument analogous to the K\H{o}v\'{a}ri-S\'{o}s-Tur\'{a}n theorem. 

In order to characterize forbidden 0-1 matrices with linear extremal functions, one may investigate an equivalent problem of determining the \textit{minimal nonlinear 0-1 matrices}, which are 0-1 matrices that have superlinear extremal functions, but all patterns strictly contained in them have extremal functions that are at most linear. Keszegh posed the problem of whether there exist infinitely many minimal nonlinear 0-1 matrices \cite{keszegh}, which was answered affirmatively by Geneson \cite{geneson2009} using a non-constructive argument. A number of properties of minimally nonlinear 0-1 matrices were identified in \cite{Crowdmath2018, GT2020}, including bounds on the ratio of their dimensions and on the number of ones in a minimally non-linear 0-1 matrix with respect to the number of rows. 

Besides identifying the minimally nonlinear 0-1 matrices, another method has been applied to the problem of characterizing the forbidden 0-1 matrices with linear extremal functions. In particular, several operations have been found which can be applied to forbidden 0-1 matrices with linear extremal functions to produce new forbidden 0-1 matrices which still have linear extremal functions \cite{tardos, keszegh}. In other words, these operations preserve linearity of the extremal function. One such linearity-preserving operation is adding rows or columns of all zeroes. In this paper, we investigate the effect of operations on saturation functions, including adding empty columns. We also note that a number of the aforementioned results have been extended to multidimensional 0-1 matrices (see, e.g., \cite{KM, GT17, geneson19, geneson21, GTT, gt23, GHLNPW}).

\subsection{The Saturation Function}

Brualdi and Cao \cite{BC} were the first to obtain results about the saturation function. In particular, they found forbidden 0-1 matrices for which the extremal function and the saturation function coincide. 

\begin{prop}[\cite{BC}]
    Let $I_k$ be the identity matrix, and $J_3 = \begin{psmallmatrix}
        & \bullet & \\
        & & \bullet \\
        \bullet & &
    \end{psmallmatrix}$. Then for all $n$, we have $\sat(n,I_k) = \ex(n,I_k)$, and $\sat(n,J_3) = \ex(n,J_3)$.
\label{thm:idmat}
\end{prop}

While this result might seem to indicate that the saturation function would often behave similarly to the extremal function, Fulek and Keszegh \cite{FK} showed that it differs markedly. While the extremal function $\ex(n,P)$ exhibits a variety of different growth rates, Fulek and Keszegh demonstrated that the saturation function exhibits much simpler, dichotomous, asymptotic behavior. In particular, they proved the following result.

\begin{thm}[\cite{FK}]
    For any pattern $P$, either $\sat(n,P) = \Theta(n)$ or $\sat(n,P) = O(1)$. Additionally, for $m_0$, $n_0$ fixed, $\sat(m_0,n,P)$ is either $\Theta(n)$ or $O(1)$, and $\sat(m,n_0,P)$ is either $\Theta(m)$ or $O(1)$.
\label{thm:dicho}
\end{thm}

The following is also asserted as part of this result:


\begin{claim}[\cite{FK}]
    If $\sat(m_0,n,P) = O(1)$ for some $m_0$, then for each $m_1 > m_0$, also $\sat(m_1, n, P) = O(1)$.
\label{claim:hwitness_assertion}
\end{claim}

The proof in \cite{FK} of Claim~\ref{claim:hwitness_assertion} of the theorem is flawed. Their strategy relies on duplicating empty rows of an $m_0 \times n$ $P$-saturating matrix, and such empty rows may not exist, as in Example~\ref{ex:FK_issue}. However, the statement is true; we offer a corrected proof in section~\ref{sec:1dsat}.

Theorem~\ref{thm:dicho} leads to a natural problem: characterize all $P$ for which the saturation function is linear (i.e.\ $\Theta(n)$), and all $P$ for which is it bounded (i.e.\ $O(1)$). Some progress has been made on this problem, which we discuss here, but the full classification has not yet been completed.

The first contributions to the classification of patterns by saturation function were from Fulek and Keszegh. Firstly, they introduced the \textit{semisaturation function} $\ssat(n,P)$, defined as follows.

\begin{definition}[Semisaturation]
    A matrix $M$ is said to be \textit{semisaturating} for a pattern $P$ if changing any entry $0$ in $M$ to a $1$ introduces a \textit{new} copy of $P$; here $M$ is allowed to contain a copy of $P$ already, which differs from the definition of saturation. The \textit{semisaturation function} $\ssat(n,P)$ is defined as the minimum possible weight of any $n \times n$ matrix that is semisaturating for $P$.
\end{definition} Clearly, if a matrix is saturating, then it is semisaturating; thus we always have $\ssat(n,P) \leq \sat(n,P)$. Fulek and Keszegh classified the growth of the semisaturation function completely as follows.

\begin{thm}[\cite{FK}]
    For any pattern $P$, we have $\ssat(n,P) = O(1)$ if and only if all of the following hold:
    \begin{enumerate}
        \item In both the first and last columns of $P$, there is a $1$ entry which is the only $1$ entry in its row.
        \item In both the first and last rows of $P$, there is a $1$ entry which is the only $1$ entry in its column.
        \item There is a $1$ entry which is the only $1$ entry in its row and column.
    \end{enumerate}
    Otherwise, $\ssat(n,P) = \Theta(n)$.
\label{thm:ssat}
\end{thm}

As $\ssat(n, P)$ is a lower bound for $\sat(n, P)$, this result also identifies a large class of patterns with $\sat(n,P) = \Theta(n)$. Seen differently, this result implies a necessary (but not sufficient) condition on $P$ for $\sat(n,P)$ to be bounded.

Fulek and Keszegh also found another class of matrices with linear saturation function. Call a matrix $M$ \textit{decomposable} if $M$ can be written as a block matrix $M = \begin{psmallmatrix} A & \mathbf{0} \\ \mathbf{0} & B \end{psmallmatrix}$ or $M = \begin{psmallmatrix} \mathbf{0} & A \\ B & \mathbf{0} \end{psmallmatrix}$ for some matrices $A$ and $B$ that are both not all-$0$ matrices, where $\mathbf{0}$ represents an all-zeroes matrix of any size.

\begin{thm}[\cite{FK}]
    Let $P$ be any decomposable pattern. Then $\sat(n,P) = \Theta(n)$.
\label{thm:decomp}
\end{thm}

Besides identifying several classes of matrices with linear saturation functions, Fulek and Keszegh also introduced a method for demonstrating that a pattern has $\sat(n,P) = O(1)$, namely by constructing a so-called \textit{witness} matrix. This construction has since been used by Geneson in \cite{geneson21a} and Berendsohn in \cite{berendsohn21} and \cite{berendsohn23} to identify examples of patterns with bounded saturation function. We adopt the terminology of Berendsohn.

Consider a pattern $P$ with no empty rows or columns. For a matrix $W$, call a row or column $i$ of $W$ \textit{expandable} with respect to $P$ if $i$ is empty and changing any entry of $i$ to a $1$ introduces a copy of $P$. A matrix $W$ is a \textit{witness} for $P$ if $W$ does not contain $P$ and $W$ has both an expandable row and an expandable column with respect to $P$.

$W$ is an \textit{explicit witness} if $W$ has an empty row and column and $W$ is saturating for $P$. Clearly all explicit witnesses are witnesses, and any witness can be used to construct an explicit witness by adding ones until it is impossible to do so without producing a copy of $P$. A \textit{horizontal (vertical) witness} for $P$ is a matrix which does not contain $P$, and which has an expandable column (row).

We can generalize the notion of witnesses to patterns which contain empty rows and columns as follows.

\begin{definition}[Witness]
    A matrix $W$ is a \textit{horizontal (vertical) witness} for a pattern $P$ containing no more than $k$ consecutive empty columns (rows) if $W$ does not contain $P$ and $W$ has $k$ consecutive expandable columns (rows). We call $W$ a \textit{witness} for $P$ if $W$ is both a horizontal and a vertical witness.
\end{definition}

These definitions allow us to formulate the following result, which Berendsohn \cite{berendsohn21} stated for patterns without empty rows or columns, but which is also clear in the general case.

\begin{thm}[\cite{berendsohn21, FK}]
    \label{thm:witness}
    Let $P$ be any pattern. Then $\sat(n,P) = O(1)$ if and only if there exists a matrix $W$ that is a witness for $P$. Similarly, $\sat(m,n_0,P) = O(1)$ if and only if there exists a vertical witness for $P$, and $\sat(m_0,n,P) = O(1)$ if and only if there exists a horizontal witness for $P$
\end{thm}

For certain matrices, the following result allows us to restrict our attention to horizontal and vertical witnesses.

\begin{thm}[\cite{berendsohn21}]
    Let $P$ be an indecomposable pattern with no empty rows or columns, with only one $1$ entry in the last row and only one $1$ entry in the last column. Then $\sat(n,P) = O(1)$ if and only if both $\sat(m_0,n,P) = O(1)$ and $\sat(m,n_0,P) = O(1)$, or equivalently, there exist both a horizontal witness $W_H$ and a vertical witness $W_V$ for $P$.
\label{thm:vhwitness}
\end{thm}

We note that this remains a necessary, but not sufficient, condition for boundedness in the general case. Clearly any witness is also both a horizontal and vertical witness, so the forward implication holds in general; however, see Example \ref{ex:q6} for a matrix with linear saturation function despite having horizontal and vertical witnesses. In this paper, we prove several new results about saturation with a fixed dimension. 

Fulek and Keszegh were able to construct a witness for one $5 \times 5$ permutation matrix, showing that it has bounded saturation function. They conjectured that other permutation matrices with bounded saturation function exist, and asked for other examples. Their conjecture was resolved in \cite{geneson21a}, where Geneson constructed an infinite class of so-called \textit{ordinary} permutation matrices, showing that all ordinary matrices have bounded saturation function and almost all permutation matrices are ordinary. 

Geneson's result has been subsumed by \cite{berendsohn23}, which completed the classification of saturation for permutation matrices as follows. Recall from Theorem \ref{thm:decomp} that all so-called \textit{decomposable} matrices have linear saturation function. Berendsohn \cite{berendsohn23} showed that, if we restrict to permutation matrices, then the converse holds as well.

\begin{thm}[\cite{berendsohn23}]
    Let $P$ be a permutation matrix. Then $\sat(n,P) = \Theta(n)$ if and only if $P$ is decomposable.
\label{thm:permut}
\end{thm}

We also note that the notion of saturation and semisaturation has been extended to multidimensional 0-1 matrices \cite{tsai23, gt23}. Containment and avoidance for $d$-dimensional 0-1 matrices are defined the same as in the $2$-dimensional case. We say that a $d$-dimensional 0-1 matrix $A$ is $P$-saturated if $A$ avoids $P$ and any matrix $A'$ obtained from $A$ by changing a zero to a one must contain $P$. The saturation function $\sat(n, P, d)$ is defined to be the minimum possible number of ones in a $P$-saturated $d$-dimensional 0-1 matrix with all dimensions of length $n$. Note that $\sat(n, P, d)$ is undefined for all $n$ sufficiently large if $P$ has no ones, and otherwise $\sat(n, P, d)$ is defined for all $n$.

We say that a $d$-dimensional 0-1 matrix $A$ is $P$-semisaturated if any matrix $A'$ obtained from $A$ by changing a zero to a one must contain a new copy of $P$. As with the $2$-dimensional case, note that this differs from the definition of saturation since it allows the matrix $A$ to contain $P$. The semisaturation function $\ssat(n, P, d)$ is defined to be the minimum possible number of ones in a $P$-semisaturated $d$-dimensional 0-1 matrix with all dimensions of length $n$. As with $2$-dimensional saturation and semisaturation functions, one can also consider non-square matrices $A$. In particular, we define $\sat(n_1, n_2, \dots, n_d, P, d)$ and $\ssat(n_1, n_2, \dots, n_d, P, d)$ analogously. It was shown in \cite{gt23} that the semisaturation function $\ssat(n, P, d)$ is always of the form $\Theta(n^r)$ for some integer $0 \le r \le d-1$. It is an open problem to determine if the same is true for the $d$-dimensional saturation function, though the same paper showed that the $d$-dimensional saturation function must be $\Omega(n)$ if it is not $O(1)$.

In this paper, we prove some new results about saturation in multidimensional 0-1 matrices. We prove a sharp upper bound on the effect of adding a new $(d-1)$-dimensional cross section to the beginning or end a forbidden $d$-dimensional 0-1 matrix $P$ on the saturation function $\sat(n, P, d)$. In particular, we show that this bound is exact in the case when we add an all-zero $(d-1)$-dimensional cross section to the beginning or end a forbidden $d$-dimensional 0-1 matrix $P$. We use this to obtain a sharp upper bound on the maximum possible value of the saturation function $\sat(n, P, d)$. We also find the first example of a $d$-dimensional 0-1 matrix with bounded saturation function for any $d > 2$. 


\section{Fixing one dimension}\label{sec:1dsat}

Throughout this section, $\mathbf0$ in block matrix notation will refer to an all-zero matrix with arbitrary dimensions that can be inferred from context. We will use some notations around witnesses. If $P = \begin{pmatrix} p_{ij} \end{pmatrix}$ and $Q = \begin{pmatrix} q_{ij} \end{pmatrix}$ are patterns where $P$ is contained in $Q$, we have an embedding $\phi$ which takes each entry of $P$ to some $q_{ij}$ in $Q$. We also have a submatrix $P'$ of $Q$, which is equal to the image of $\phi$. If $W$ is a horizontal (vertical) witness for $P$ with a single expandable column $j$ (row $i$), then we let $W_i$ ($W_j$) denote the matrix obtained by changing the entry $w_{ij}$ to a $1$. When $P$ is contained in $W_i$, we let $P_i$ denote some (not necessarily unique) copy of $P$ in $W_i$. The embedding that takes $P$ to $P_i$ in $W_i$ is $\phi_i$. There is also a submatrix $P_i'$ of $W$, which contains the same rows and columns as $P_i$ in $W_i$. The submatrix $P_i'$ differs from $P_i$ only by a single entry $o'$ in row $i$ of the expandable column: the entry $o'$ is a $1$ in $P_i$, but a $0$ in $P_i'$. This entry corresponds to some one entry $o$ of $P$ itself; that is, $o' = \phi_i(o)$; define $_{(i)}P$ by changing $o$ to a zero. Then the image of $_{(i)}P$ under $\phi_i$ lies inside of $P_i$, but with the one entry in the expandable column removed; thus all of the one entries of $_{(i)}P$ correspond under $\phi_i$ to one entries which are in fact in the original witness $W$. With slight abuse of notation, we also use $\phi_i$ to represent this embedding which takes $_{(i)}P$ into $W$ itself (understanding that this is properly a restriction of $\phi_i$ as first defined to $_{(i)}P$).

We first develop the simpler case of $\sat(m_0,n,P)$. For all of the results in this section, we omit symmetric phrasing for $\sat(m,n_0,P)$, for the sake of brevity. However, all of these results apply similarly to this case.

Theorem \ref{thm:ssat} gives a necessary and sufficient condition for $\ssat(n,P)$ to be bounded. We state a similar condition characterizing when $\ssat(m_0,n,P) = O(1)$.
\begin{prop}
\label{prop:vhssat}
    For a pattern $P$ and sufficiently large $m_0$, the semisaturation function $\ssat(m_0, n, P) = O(1)$ if and only if the top and bottom rows of $P$ each contain a $1$ entry which is the only $1$ entry in its column.
    %
\end{prop}
\begin{proof}
    Let $P$ have dimensions $k\times\ell$. Suppose that the first and last rows of $P$ both contain a $1$ entry that is the only $1$ in its column. We will call these $1$ entries \textit{special} $1$ entries. Then for $m_0\ge2(k-1)$, consider the $m_0 \times n$ matrix $M$ such that its first and last $\ell-1$ columns are all $1$s, with $0$s elsewhere. That is, \[M=\begin{pmatrix}\mathbf1_{m_0\times(\ell-1)}&\mathbf0&\mathbf1_{m_0\times(\ell-1)}\end{pmatrix}\]
    where $\mathbf1_{m_0\times(\ell-1)}$ is the $m_0\times(\ell-1)$ matrix of all $1$s.
    For this matrix, if we change a $0$ that is not in the last $k-1$ rows into a $1$, then the new $1$ entry plays the role of the special $1$ entry in the first row of $P$ in a new copy of the pattern in $M$. Similarly, if the $0$ is not in the first $k-1$ rows of $M$, then the new $1$ entry assumes the role of the special $1$ entry in the last row of $P$ in some new copy of $P$ in $M$. Since the weight of $M$ is $2m_0(\ell-1)$, $\ssat(m_0,n,P)\le2m_0(\ell-1)=O(1)$.

    Now suppose that $P$ does not have a special $1$ entry in either its first or last row, and let $M$ be an $m_0 \times n$ matrix that is semisaturated for $P$. Then, every column of $M$ would have to contain at least one $1$ entry. If there was an empty column, then changing a $0$ in that column into a $1$ in the first row would not create a new copy of $P$ because the corresponding $1$ entry in $P$ would have to be special, which is impossible. Therefore, there are at least $n$ $1$ entries in $M$, so $\ssat(m_0,n,P)=\Omega(n)$. Since $\ssat(m_0,n,P)\le\sat(m_0,n,P)=O(n)$, we conclude $\ssat(m_0,n,P)=\Theta(n)$.
%
\end{proof}

\subsection{A proof for bounded saturation functions}

The statement of Theorem \ref{thm:dicho} as originally published in \cite{FK} contained the following assertion. We note that the proof given in \cite{FK} fails to fully cover this statement, and provide a corrected proof.

\begin{thm}[\cite{FK}]
    If $\sat(m_0,n,P) = O(1)$ for some $m_0$, then for each $m_1 > m_0$, also $\sat(m_1, n, P) = O(1)$.
\label{thm:hwitness_correction}
\end{thm}

The proof given in \cite{FK} of this fact proceeds by finding empty rows in a $P$-saturating $m_0 \times n$ matrix and multiplying these empty rows to find a $P$-saturating $m_1 \times n$ matrix of the same weight; it follows that if the former weight is bounded, the latter must also be. However, this strategy relies on the original $P$-saturating matrix having an empty row. This cannot work in general.
\begin{example}\label{ex:FK_issue}
    The pattern \[P = \begin{pmatrix}
    \bullet & \bullet & & \bullet \\ \bullet & & \bullet & \bullet
    \end{pmatrix}\] has no $1$ entries which are alone in their row; thus no horizontal witness for $P$ can possibly have an empty row. Thus the strategy of Fulek and Keszegh, which relies on empty rows, would certainly fail on $P$; but by Proposition \ref{prop:intermediaries_1}, we will see that $P$ does have a horizontal witness.
\end{example}

We give a corrected proof of Theorem \ref{thm:hwitness_correction}, using the similar idea of multiplying the rows of $P$.

\begin{proof}[Proof of Theorem~\ref{thm:hwitness_correction}]
    Suppose for some $k \times l$ pattern $P$, the function $\sat(m_0,n,P) = O(1)$. Then for some $n$, we have an $m_0 \times n$ horizontal witness $W = (w_{ij})$ for $P$, with an expandable column $j$. In particular, when we change the entry $w_{m_0j}$ to a $1$, we introduce a copy $P_{m_0}$ of $P$ in $W$. Let the number of $1$ entries in the bottom row of $P$ be $b$. Then the bottom row of $P_{m_0}$ contains $1$ entries $w_{m_0j_1},\dots,w_{m_0j_a},\dots,w_{m_0j_b}$, where $j_a = j$ is the expandable column.
    
    Define the matrix $W'$ by appending to $W$ a row with $b - 1$ one entries located in the columns $j_1, \dots, \widehat{j_a}, \dots, j_b$ (the notation meaning that the column $j_a = j$ is omitted). If $W'$ contains a copy of $P$, then since $W$ avoids $P$, this copy must contain the bottom row of $W'$, but this row has only $b - 1 < b$ one entries, so certainly $W'$ avoids $P$.  Furthermore, if we change the entry of $W'$ in the bottom row and column $j$ to a $1$, we obtain a copy of $P$ as follows: considering $P_{m_0}$ as above, replace the bottom row of $P_{m_0}$ with the bottom row of $W'$. By construction this has $1$ entries in all the same columns, and so it is a copy. Then column $j$ remains expandable.

    We have seen that $W'$ is an $(m_0 + 1) \times n$ horizontal witness for $P$. Repeating this argument inductively gives an $m_1 \times n$ horizontal witness for any integer $m_1 > m_0$, which is enough to imply the claim.
\end{proof}

Fulek and Keszegh asked in \cite{FK} whether a finite algorithm exists to determine whether a pattern $P$ has bounded $\sat(n,P)$. As they note, this is equivalent to asking whether there exists a computable function $f(k,l)$ such that, for a $k \times l$ pattern $P$, there exists a witness for $P$ of dimensions at most $f(k,l) \times f(k,l)$. We resolve a version of this question for $\sat(m_0,n,P)$.

\begin{prop}
    Let $P$ be any $k \times l$ pattern with no empty columns. If $\sat(m_0,n,P) = O(1)$, then there is a horizontal witness for $P$ with at most $(l-1)m_0 + 1$ columns, and $\sat(m_0,n,P) \leq (l-1)m_0^2$.
\end{prop}
\begin{proof}
    That the latter claim follows from the former is clear; we show the former. We know we have some horizontal witness $W$ for $P$, with an expandable column. Clearly, if a column $j$ is contained in none of the patterns $P_i$ in $W$, then $j$ may be deleted and $W$ will remain a witness. We can therefore assume that every column in $W$ is in some $P_i$. The bound follows by seeing that there are $m_0$ patterns $P_i$, each of which has $l$ columns, and one of these columns is the expandable column, which is the same for every $P_i$.
\end{proof}

This bound is in terms of $m_0$; however, we fail to find a bound independent of $m_0$. Such a function would resolve the question of Fulek and Keszegh. Since, for most patterns discussed here and in Section \ref{sec:patterns}, we have either found a fairly small witness or shown that none exists, it seems likely that such a bound may exist.

\subsection{A generalization of Theorem \ref{thm:vhwitness}}

Theorem \ref{thm:vhwitness} allows us to glue together a horizontal and vertical witness for certain matrices, reducing the study of saturation for square matrices to the case with one dimension fixed. We extend the result in Theorem \ref{thm:vhwitness} to apply to a class of patterns that do not have only one $1$ entry in their last row and column. 
\begin{definition}[Strong Indecomposability]
    We define an indecomposable pattern $P$ to be strongly indecomposable if there do not exist matrices $P_1$, $P_2$, $P_3$, and $P_4$, at most one of which is empty, for which $P$ can be represented as \begin{align*}
        P = &\begin{pmatrix}
            \mathbf 0 & P_1 & \mathbf 0 \\
            P_2 & \mathbf 0 & P_3 \\
            \mathbf 0 & P_4 & \mathbf 0
        \end{pmatrix}.
    \end{align*}
\end{definition}

We present a stronger necessary and sufficient condition for strongly indecomposable matrices to have bounded saturation function. 
\begin{thm}\label{thm:str_indecomp}
    For a given strongly indecomposable pattern $P$, $\sat(n,P) = O(1)$ if and only if there exists a vertical witness $W_V$ and a horizontal witness $W_H$ corresponding to $P$, each such that at least one zero in the expandable row (column) corresponds to a $1$ in $P$ that is alone in its column (row). 
\end{thm}
\begin{proof}
    Firstly, we show that this condition is necessary for a bounded saturation function. If there does not exist a vertical witness $W_V$ such that at least one zero in the expandable row corresponds to a $1$ in $P$ that is alone in its column, then every zero in an empty row in a saturating matrix for $P$ must correspond to a $1$ in $P$ that is not alone in its column. Thus for every zero in the empty row, there must exist another $1$ in its column. Thus no matrix that is saturating for $P$ can contain both an empty row and an empty column, and the function $\sat(n,P)$ cannot be bounded as $n$ grows. A similar argument holds for $W_H$.

    Second, we show that this condition is sufficient for a bounded saturation function by constructing a witness $W$ for $P$ given both a vertical witness $W_V$ and a horizontal witness $W_H$ with the property outlined. For $W_H$, if row $n$ contains a zero which corresponds to a $1$ in $P$ alone in its row, then we can take the matrix $A$ to be the submatrix of $W_H$ consisting of the first $n$ rows, and matrix $D$ to be the submatrix of $W_H$ consisting of the remaining rows. Similarly, if column $m$ of $W_V$ contains a zero corresponding to a $1$ in $P$ alone in its column, then we define $B$ to be the submatrix of $W_V$ consisting of the first $m$ columns, and $C$ the submatrix of $W_V$ consisting of the remaining rows. 

    Let \begin{align*}
        W = &\begin{pmatrix}
            \mathbf 0 & A & \mathbf 0 \\
            B & \mathbf 0 & C \\
            \mathbf 0 & D & \mathbf 0
        \end{pmatrix}.
    \end{align*} We seek to prove that $W$ is a witness for $P$. This is equivalent to the following two claims: That $W$ contains both an expandable empty row and an expandable empty column, and that $W$ does not contain $P$ as a submatrix. 

    It is trivially true that $W$ contains both an expandable empty row and an expandable empty column. Appending any number of empty rows to $W_H$ directly below row $n$ will preserve the condition of the expandable empty column because all of the zeroes that already existed in the empty column will still have the property that changing them into a one will create a copy of $P$. Meanwhile, all of the zeroes created immediately below row $n$ will form copies of $P$ with the same ones as the zero in row $n$. Similar logic applies for $W_V$'s expandable empty row. 
    
    Appending columns, empty or not, cannot transform a matrix with an expandable empty column into a matrix without one, and appending rows cannot transform a matrix with an expandable empty row into a matrix without one, so the process of transforming $W_H$ into $W$, first by appending empty rows underneath row $n$ and then by appending columns to the left and right, will not break the condition of there being an expandable empty column. Similarly, the process of transforming $W_V$ into $W$ will not break the condition of there being an expandable empty row.

    It therefore suffices to show that $W$ does not contain any copies of $P$. For this, we proceed with a proof by contradiction. 

    Assume that $P$ is contained within $W$. Any nontrivial pattern $W'$ that is contained within $W$ must be of one of the following forms:

    \begin{enumerate}
        \item $W'$ contains ones in only one of $\{A, B, C, D\}$.
        \item $W'$ contains ones in $A$ and $D$ or $B$ and $C$. 
        \item $W'$ contains ones in $A$ and $B$, $A$ and $C$, $B$ and $D$, or $C$ and $D$.
        \item $W'$ contains ones in any three of $\{A, B, C, D\}$.
        \item $W'$ contains ones in all of $\{A, B, C, D\}$.
    \end{enumerate}

    If $W'$ is of the first or second form, then $W'$ cannot be $P$ because, by definition, $W_V$ and $W_H$ do not contain $P$. 
    
    If $W'$ is of the third form, then we assume that it contains ones in $A$ and $B$, then proceed with identical reasoning for the other three cases. Its ones can be split into two subsets, $S_1$ and $S_2$, with $S_1$ corresponding to the ones in $A$ and $S_2$ corresponding to the ones in $B$. Because every row in $A$ is above every row in $B$, and because every column in $A$ is to the right of every column of $B$, every one in $S_1$ is above and to the right of every one in $S_2$, and $W'$ is not indecomposable by definition. Since $P$ is indecomposable, $W' \neq P$.

    If $W'$ is of the fourth form, for example if its ones belong to $\{A, B, C\}$, we, like in the previous case, split its ones into three subsets, $S_1$, $S_2$, and $S_3$, belonging to $A$, $B$, and $C$, respectively. All of the ones in $S_1$ are above all of the ones in $S_2$ and $S_3$. All of the ones in $S_2$ are to the left of $S_1$ and $S_3$. All of the ones in $S_3$ are to the right of all of the ones in $S_1$ and $S_2$. Therefore, $W'$ is of the form  \begin{align*}
        &\begin{pmatrix}
            \mathbf{0} & P_1 & \mathbf{0} \\
            P_2 & \mathbf{0} & P_3 \\
        \end{pmatrix}.
    \end{align*} Since $P$ is strongly indecomposable, $W' \neq P$. 

    Finally, if $W'$ is of the fifth form, by the same logic, $W'$ is of the form \begin{align*}
        &\begin{pmatrix}
            \mathbf{0} & P_1 & \mathbf{0} \\
            P_2 & \mathbf{0} & P_3 \\
            \mathbf{0} & P_4 & \mathbf{0}\\
        \end{pmatrix},
    \end{align*} and $W' \neq P$.

    Thus $P$ is not contained within $W$, and $W$ is a witness for $P$. 
    
    Therefore, this condition is both necessary and sufficient for matrices of this form. 
\end{proof}
\label{sec:operations}

\section{Inserting rows and columns}
\label{sec:insert}

In this section, we show several results on the effect of inserting columns and rows into a pattern. In the case when the columns and rows are empty, we obtain an exact result.

\begin{thm}\label{thm_col_add_op}
    Let $P$ be any pattern, and let $P'$ be any pattern derived from $P$ by adding a single column before the first column of $P$. Then $\sat(m,n,P') \le m + \sat(m,n-1,P)$.
\end{thm}

\begin{proof}
    Let $M$ be an $m \times (n-1)$ matrix that is saturating for $P$. Consider the matrix $M'$ derived from $M$ by adding a column of all $1$ entries before $M$. We claim that $M'$ is saturating for $P'$. First, we prove by contradiction that $M'$ avoids $P'$. Suppose for contradiction that $M'$ contains $P'$. Removing the leftmost column of $P'$, we obtain a copy of $P$ in $M'$; furthermore this copy of $P$ does not include the first column of $M'$, so $M$ also contains a copy of $P$, which is a contradiction. 
    
    To see that $M'$ is saturating, we now show that adding any $1$ to $M'$ introduces a copy of $P'$. Since the first column of $M'$ is all $1$, any $1$ entry must be added outside the first column. Then this $1$ entry is in $M$, so it is part of a copy of $P$ in $M$; this copy of $P$ together with the first column of $M'$ forms a copy of $P'$ in $M'$ as desired. Thus $\sat(m,n,P') \leq m + \sat(m,n-1,P)$.
\end{proof}

\begin{cor}\label{cor_empty_col_add_op}
    Let $P$ be any pattern, and let the pattern $P'$ be derived from $P$ by adding an empty column before the first column of $P$. Then $\sat(m,n,P') = m + \sat(m,n-1,P)$.
\end{cor}

\begin{proof}
By Theorem~\ref{thm_col_add_op}, we have $\sat(m,n,P') \leq m + \sat(m,n-1,P)$. For the other direction, now let $M$ be any $m \times n$ matrix saturating for $P'$. If the first column of $M$ has any $0$ entries, then changing them to a $1$ cannot introduce a copy of $P'$, since the first column of $P'$ is empty. Thus the first column of $M$ must have all ones. We claim that the matrix $M'$ derived from $M$ by deleting the first column is saturating for $P$. Clearly if $M'$ contains a copy of $P$, then this copy of $P$ with the first column of $M$ added forms a copy of $P'$ in $M$, so $M'$ cannot contain a copy of $P$. If a $0$ entry in $M'$ is changed to a $1$, we know that this creates a copy of $P'$ in $M$; then this also introduces a copy of $P$ in $M$ which cannot include the first column, which is a copy of $P$ in $M'$. Thus, $M'$ is saturating for $P$, so $\sat(m,n-1,P) \leq \sat(m,n,P') - m$. Combining the two directions of inequality, the result follows.
\end{proof}

\begin{cor}
    There exist patterns $P, P'$ with $\sat(n,P) = O(1)$ and $\sat(n,P') = \Theta(n)$ for which $P'$ is obtained from $P$ by adding an empty column.
\end{cor}

\begin{proof}
    Let $P$ be any pattern with $\sat(n,P) = O(1)$, e.g., almost any permutation matrix, and let $P'$ be obtained from $P$ by adding an empty column on the left. Then, $\sat(n,P') = \Theta(n)$ by Corollary~\ref{cor_empty_col_add_op}.
\end{proof}

We begin with an easy statement, and in Propositions \ref{prop:empty_col_h} and \ref{prop:emptyrows} we investigate its converse.
\begin{prop}
    Let $P$ be a pattern and $P'$ be obtained from $P$ by adding empty rows and columns. If $\sat(m,n,P')=O(1)$, then $\sat(m,n,P)=O(1)$ as well.
\end{prop}
\begin{proof}
    Let $W'$ be a witness for $P'$. Since $W'$ does not contain $P'$, it does not contain $P$ either. Since changing each entry in an expandable column of $W'$ induces a copy of $P'$ in the resulting matrix and $P'$ contains $P$, we also obtain a copy of $P$ in the resulting matrix. Therefore, $W'$ is also a witness for $P$ and $\sat(m_0,n,P)=O(1)$.
\end{proof}

It remains unknown whether there exist patterns $P$, $P'$ with $\sat(n, P) = O(1)$ and $\sat(n, P') = \Theta(n)$, for which $P'$ is derived from $P$ by adding an empty \textit{interior} column, i.e., between the first and last columns of $P$. We can, however, show the following for $\sat(m_0, n, P)$.

\begin{prop}
\label{prop:empty_col_h}
    Let $P$ be a pattern without empty columns, and let $P'$ be derived from $P$ by adding empty columns. Then $\sat(m_0, n, P') = O(1)$ if $\sat(m_0, n, P) = O(1)$.
\end{prop}

\begin{proof}
    Assume first that $\sat(m_0,n,P)=O(1)$. From Theorem \ref{thm:witness} we know there exists a horizontal witness $W = \begin{pmatrix} a_{ij} \end{pmatrix}$ for $P$; let $W$ be $m \times n$. It suffices to give a horizontal witness for $P'$. Suppose $P'$ contains up to $k$ consecutive empty columns. Construct the $m \times ((2k+1)(n-1) + 1)$ matrix $W' = \begin{pmatrix} b_{ij} \end{pmatrix}$ by inserting $2k$ empty columns between each of the columns of $W$; that is, let
    \[b_{ij'} = \begin{cases}
        a_{ij} & j' = 1 + (2k+1)(j - 1) \\
        0 & \text{else}
    \end{cases}.\]
    Clearly $W'$ does not contain a copy of $P'$, as otherwise $W$ would contain $P$. Let the $j$th column of $W$ be expandable. Let $j'=(2k+1)(j-1)+1$. We firstly see that the $j'$th column of $W'$ is expandable; if we flip $b_{ij'}$ to a $1$, then because $W$ is a witness, we introduce a copy of $P$ which is contained entirely within the nonempty columns of $W'$. Then as there are $2k$ empty columns between each nonempty column of $W'$, we are able to introduce enough empty columns into this copy of $P$ to find a copy of $P'$ in $W'$.
    
    We can now show further that the $j' - k$th through $j' + k$th columns of $W'$ are all expandable. Clearly these are all empty. Consider the matrix $W''$ derived by adding a $1$ at $b_{ij''}$, where $j'' \in [j' - k, j' + k]$. We know the $j'$th column is expandable, so in particular, when we flip $b_{ij'}$ to a $1$, we introduce a copy $P'_0$ of $P'$.

    We modify $P'_0$ to create a copy of $P'$ with a 1 entry at $b_{ij''}$ rather than $b_{ij'}$. Clearly the entry $b_{ij'}$ is the only $1$ entry in its column, as is $b_{ij''}$, and they lie in the same row. Furthermore, we have only at most $k$ empty columns in $P'_0$ between column $j'$ and another nonempty column, and we know column $j''$ also has at least $k$ empty columns between itself and another nonempty column. Suppose there are $k_1$ empty columns to the left of $j'$ in $P'$, and $k_2$ to the right; then we can create $P'_1$ by taking the same columns and rows as $P'_0$, except that column $j'$ is replaced by $j''$, and these $k_1$ and $k_2$ empty columns are replaced by $k_1$ empty columns immediately to the left and $k_2$ empty columns immediately to the right of column $j''$. We observe that $P'_1$ is a copy of $P'$ in $W'$. Thus we may safely construct a copy of $P'$ including $b_{ij''}$. 
    
    Therefore in particular the column $j''$ is expandable. Then $W'$ contains $2k + 1 > k$ consecutive expandable columns, confirming the other condition for $W'$ to be a horizontal witness for $P'$. Thus $\sat(m_0,n,P') = O(1)$.
\end{proof}

The corresponding statement of Proposition \ref{prop:empty_col_h} for adding empty rows is not as clear in general; we require certain assumptions on the witness for $P$.

\begin{prop}
\label{prop:emptyrows}
    Let $P$ be a pattern without empty rows and $\sat(m_0,n,P)=O(1)$, and let $W$ be a horizontal witness for $P$ containing an expandable column $j$ such that each $w_{ij}$ is the only $1$ entry in its row of the submatrix $P_i$ of $W_i$. Let $P'$ be derived from $P$ by adding empty rows in the interior (i.e. strictly between the first and last row). Then $\sat(m_0,n,P') = O(1)$.
\end{prop}

\begin{proof}
    Suppose that $P'$ is derived by adding up to $k$ consecutive empty rows. Let $W'$ be derived from $W$ by adding $2k - 1$ consecutive empty rows between each row of $W$. We have an embedding $\psi$ of $W$ into $W'$. It is clear that if $W'$ contains $P'$, then $W$ contains $P$. We claim that column $j$ is still expandable. Consider any row $i'$ of $W'$. Let $\psi(i)$ be the closest row of $W'$ to $i'$ which corresponds to a row of $W$; without loss of generality suppose that $\psi(i) \leq i'$. Then $\psi(P_i)$ can be extended to a copy $P'_{\psi(i)}$ of $P'$ by adding to it empty rows. This copy $P'_{\psi(i)}$ has at most $k$ empty rows consecutively after the row $\psi(i)$. We know $i' - \psi(i) \leq k - 1$, so there are at least $k$ empty rows in $W'$ between $i'$ and $\psi(i+1)$, the next row corresponding to a row of $W$. Thus, we may build $P'_{i'}$ from $P'_{\psi(i)}$ by replacing the row $\psi(i)$ with $i'$, and the up to $k$ consecutive empty rows after $\psi(i)$ with the consecutive empty rows immediately after $i'$. This gives a copy $P'_{i'}$ in $W'_{i'}$, so the column $j$ is expandable, so $W'$ is a horizontal witness for $P'$.
\end{proof}

We can combine these results to see the following.

\begin{cor}
    Let $P$ be a pattern with bounded saturation function satisfying the hypotheses of $\ref{thm:vhwitness}$ such that each $1$ entry $p_{ij}$ is the only $1$ entry in row $i$ if and only if it is the only $1$ entry in column $j$. Let the pattern $P'$ be derived from $P$ by adding empty rows and columns in the interior (i.e., strictly between the first and last row and column). Then $\sat(n,P') = O(1)$.
\end{cor}

\begin{proof}
    Because $P$ satisfies the hypotheses of $\ref{thm:vhwitness}$, also $P'$ must; thus it suffices to give horizontal and vertical witnesses for $P'$. We know there exists a horizontal witness for $P$; furthermore, by our condition on $P$, this witness must satisfy the hypotheses of $\ref{prop:emptyrows}$. Thus, if $P''$ is derived from $P$ by adding empty rows, we can find a horizontal witness for $P''$. Then by $\ref{prop:empty_col_h}$ we can find a horizontal witness for $P'$. Similarly, there exists a vertical witness for $P'$, which completes the proof.
\end{proof}

The above result covers in particular all permutation matrices $P$. More generally, we ask whether the hypotheses are strictly necessary; in particular, we make the following conjecture, which would widely generalize the results of this section.

\begin{conj}
\label{conj:empty_rows_cols}
    Let $P$ be any pattern, and let $P'$ be derived from $P$ by adding empty rows and columns in the interior. Then, $\sat(n,P') = O(\sat(n,P))$ and $\sat(m_0,n,P') = O(\sat(m_0,n,P))$.
\end{conj}

This would reduce the study of the asymptotic behavior of the saturation function to the case with no empty interior rows or columns.

\section{The Kronecker Product}
\label{sec:kron}

For two matrices $P = \begin{pmatrix} p_{i,j} \end{pmatrix}$, $Q = \begin{pmatrix} q_{i,j} \end{pmatrix}$, with dimensions $k \times l$, $m \times n$ respectively, we consider the saturation of the $km \times ln$ Kronecker product $P \otimes Q = \begin{pmatrix} r_{i,j} \end{pmatrix}$, defined by \[
    r_{i,j} = p_{a,c}q_{b,d},
\] where \begin{align*}
    i &= m(a-1) + b \\
    j &= n(c-1) + d.
\end{align*}

We show the following relationship for the semisaturation function.

\begin{thm}\label{thm:kron_prod}
    Let $P = \begin{pmatrix} p_{i,j} \end{pmatrix}$, $Q = \begin{pmatrix} q_{i,j} \end{pmatrix}$ be any patterns, where $P$ is $k \times l$, $Q$ is $m \times n$. Then the semisaturation function $\ssat(n, P \otimes Q) = O(1)$ if and only if both $\ssat(n, P) = O(1)$ and $\ssat(n, Q) = O(1)$. Similarly, $\ssat(m_0, n, P \otimes Q) = O(1)$ if and only if both $\ssat(m_0, n, P) = O(1)$ and $\ssat(m_0, n, Q) = O(1)$.
\end{thm}

\begin{proof}
    For the forward direction, suppose that $\ssat(n, P \otimes Q) = O(1)$. Then by Theorem \ref{thm:ssat}, there exists a $1$ entry $r_{i,1}$ in the first column of $P \otimes Q$ which is the only $1$ entry in its row. Write $i = m(a-1) + b$; then by definition of the Kronecker product, $r_{i,1} = p_{a,1}q_{b,1}$. Thus we must have $p_{a,1} = q_{b,1} = 1$. Suppose for some $c > 1$ that we have $p_{a,c} = 1$. Then letting $j = n(c-1) + 1$, we have $r_{i,j} = p_{a,c}q_{b,1} = 1$, which contradicts the fact that $r_{i,1}$ was the only $1$ entry in its row in $P \otimes Q$. Thus $p_{a,1}$ is the only $1$ entry in its row in $P$. Similarly we see that $q_{b,1}$ is the only $1$ entry in its row in $Q$. Thus both $P$ and $Q$ have a $1$ entry in the first column which is the only $1$ entry in its first row. Similar reasoning shows that both $P$ and $Q$ have
    \begin{enumerate}[label=(\roman*)]
        \item a $1$ entry in the last column which is the only entry in its row,
        \item a $1$ entry in the first row which is the only entry in its column, and
        \item a $1$ entry in the last row which is the only entry in its column.
    \end{enumerate}

    To show $P$ and $Q$ have bounded semisaturation function, by Theorem \ref{thm:ssat}, it remains to show that both contain a $1$ entry which is the only $1$ in its row and column. We know that $P \otimes Q$ contains such a $1$ entry $r_{i,j}$. We write $i = m(a-1) + b$, $j = n(c-1) + d$. Then $r_{i,j} = p_{a,c}q_{b,d}$. Then $p_{a,c}=q_{b,d}=1$. Suppose for some $a' \neq a$ that the entry $p_{a',c} = 1$. Then letting $i' = m(a' - 1) + b \neq i$, we see $r_{i',j} = 1$, contradicting the fact that $r_{i,j}$ is the only $1$ entry in its column. Thus $p_{a,c}$ is the only $1$ entry its column (and similarly in its row) in $P$. Similarly also $q_{b,d}$ is the only $1$ entry in its column and row in $Q$.

    Thus $P$ and $Q$ both have bounded semisaturation function, completing the forward direction of the proof.

    For the reverse implication, suppose $P$ and $Q$ both have bounded semisaturation function. Then we have $1$ entries $p_{a,1}$, $q_{b,1}$ in $P$ and $Q$ respectively which are the only $1$ entries in rows $a$ and $b$ respectively. Taking $i = m(a-1) + b$, it is easy to see that the entry $r_{i,1}$ in $P \otimes Q$ is a $1$ entry which is the only $1$ in its row; we can symmetrically check the other conditions of Theorem \ref{thm:ssat}, except for the final condition, which we check separately. We have $1$ entries $p_{a,c}$, $q_{b,d}$ which are the only $1$ entries in their respective rows and columns in $P$ and $Q$; then taking $i = m(a-1) + b$, $j = n(c-1) + d$, it is again easy to see that the entry $r_{i,j}$ in $P \otimes Q$ is the only $1$ entry in its row and column. Thus $\ssat(n, P \otimes Q) = O(1)$.

    We can use very similar reasoning by Proposition \ref{prop:vhssat} to prove the version for $\ssat(m_0, n, \cdot)$.
\end{proof}

We have the following immediate corollary.


\begin{cor}{\label{cor:kron_prod_flip}}
    For any patterns $P$ and $Q$, $\ssat(n, P \otimes Q) = O(1)$ if and only if $\ssat(n, Q \otimes P) = O(1)$.
\end{cor}

\begin{proof}
    Suppose that $\ssat(n, P \otimes Q) = O(1)$. Then by Theorem \ref{thm:kron_prod}, we have $\ssat(n, P) = O(1)$ and $\ssat(n, Q) = O(1)$, so again by Theorem \ref{thm:kron_prod}, we have $\ssat(n, Q \otimes P) = O(1)$.
\end{proof}

We note that both Theorem \ref{thm:kron_prod} and Corollary \ref{cor:kron_prod_flip} fail when analogous statements are made for saturation, as illustrated by the following example.

\begin{example}
    Let \[P = \begin{psmallmatrix}
        & \bullet & & \\ & & & \bullet \\ \bullet & & & \\ & & \bullet &
    \end{psmallmatrix}, \: Q = \begin{psmallmatrix} \bullet & \\ & \bullet \end{psmallmatrix}.\] Then $\sat(n,P \otimes Q) = O(1)$, while $\sat(n,Q) = \Theta(n)$. Furthermore, $\sat(n, Q \otimes P) = \Theta(n)$.
\end{example}

All of the above follow from Theorem \ref{thm:permut}. This demonstrates the failure of Corollary \ref{cor:kron_prod_flip} and of the forward direction of Theorem \ref{thm:kron_prod} when analogous statements are made for saturation. However, we fail to resolve the reverse implication of Theorem \ref{thm:kron_prod} when analogous statements are made for saturation. Based on the known case for permutation matrices, we formulate the following conjecture.

\begin{conj}
\label{conj:kron}
    If $P$ and $Q$ are patterns with $\sat(n,P) = O(1)$ and $\ssat(n,Q) = O(1)$, then $\sat(n, P \otimes Q) = O(1)$.
\end{conj}

A special case of Proposition \ref{prop:empty_col_h} comes from considering the Kronecker product $P' = P \otimes \begin{pmatrix}
    1&0&\cdots&0
\end{pmatrix}$, or more generally any Kronecker product of $P$ with a $1 \times k$ matrix containing exactly one $1$ entry. This operation is equivalent to adding empty columns to $P$. This yields the following observation.

\begin{obs}
    Proposition \ref{prop:empty_col_h} implies a version of Conjecture \ref{conj:kron} for $\sat(m_0,n,\cdot)$ in the case of a $1 \times k$ matrix $Q$ with exactly one $1$ entry.
\end{obs}

\section{Saturation in Multidimensional Matrices}\label{sec:ddim}

In this section we adopt the terminology of \cite{tsai23}, alongside the following conventions. Fix a $d$-dimensional $n_1 \times \cdots \times n_d$ pattern $P$. A \textit{layer} is a $(d-1)$-dimensional cross section of $P$ which fixes the final coordinate $x_d = k$ for some $k \in \{1,\dots,n_d\}$. We also refer specifically to \textit{layer} $k$, denoted $P_k$. More generally, we define an $i$-layer of $P$ to be a $(d-1)$-dimensional cross section of $P$ which fixes the $i^{\text{th}}$ coordinate $x_i = k$ for some $k \in \{1,\dots,n_i\}$.

Let $A$ be the $4 \times 4 \times 6$ pattern whose layers are given by the following matrices:

\begin{align*}
    A_1 &= \begin{pmatrix}
        & & \bullet & \\ & & & \bullet \\ \bullet & & & \\ & \bullet & &
    \end{pmatrix} &
    A_2 &= \begin{pmatrix}
        & & & \\ \textcolor{white}{\bullet} & \textcolor{white}{\bullet} & \textcolor{white}{\bullet} & \textcolor{white}{\bullet} \\ & & & \\ & & & \bullet \\
    \end{pmatrix} &
    A_3 &= \begin{pmatrix}
        \bullet & & & \\ \textcolor{white}{\bullet} & \textcolor{white}{\bullet} & \textcolor{white}{\bullet} & \textcolor{white}{\bullet} \\ & & & \\ & & & \\
    \end{pmatrix} \\
    A_4 &= \begin{pmatrix}
        & & & \\ \textcolor{white}{\bullet} & \textcolor{white}{\bullet} & \textcolor{white}{\bullet} & \textcolor{white}{\bullet} \\ & & & \\ \bullet & & & \\
    \end{pmatrix} &
    A_5 &= \begin{pmatrix}
        & & & \bullet \\ \textcolor{white}{\bullet} & \textcolor{white}{\bullet} & \textcolor{white}{\bullet} & \textcolor{white}{\bullet} \\ & & & \\ & & & \\
    \end{pmatrix} &
    A_6 &= \begin{pmatrix}
        & \bullet & & \\ \bullet & & & \\ & & & \bullet \\ & & \bullet &
    \end{pmatrix}.
\end{align*}

Note that $A$ has no one entry alone in every cross-section. Thus, we have $\sat(n,A) = \Theta(n)$ according to the criteria in \cite{tsai23}. However, if we fix the first two dimensions to have length $6$ and we let the number of layers vary, we show that the resulting saturation function is bounded.

\begin{theorem} \label{thm:3d_bounded_sat}
    The pattern $A$ has saturation function $\sat(6,6,n;A) = O(1)$.
\end{theorem}

\begin{proof} For any $n \geq 8$, consider the $6 \times 6 \times n$ pattern $W$, with layers as follows:

\begin{align*}
    W_1 &= \begin{pmatrix}
        & & \bullet & & & \\
        & & \bullet & & & \\
        & & \bullet & & & \\
        & & & \bullet & \bullet & \bullet \\
        \bullet & & & & & \\
        & \bullet & & & &
    \end{pmatrix} &
    W_2 &= \begin{pmatrix}
        & & & & \bullet & \\
        & & & & & \bullet \\
        \bullet & \bullet & \bullet & & & \\
        & & & \bullet & & \\
        & & & \bullet & & \\
        & & & \bullet & \bullet & \bullet
    \end{pmatrix} \\
    W_3 &= \begin{pmatrix}
        \bullet & & & & \bullet & \\
        \bullet & & & & \bullet & \\
        \bullet & & & & \bullet & \\
        & & & & & \bullet\\
        \bullet & \bullet & \bullet & & & \bullet\\
        & & & \bullet & & \bullet
    \end{pmatrix} &
    W_4 &= \begin{pmatrix}
        \bullet & \bullet & \bullet & & & \\
        & & & \bullet & \bullet & \bullet \\
        \bullet & & & & & \\
        & \bullet & & & & \\
        & \bullet & & & & \\
        \bullet & \bullet & & & & \bullet
    \end{pmatrix}.
\end{align*}

For $k \in \{1,2,3,4\}$, we let $W_{n-k+1}$ be the reflection of $W_k$ across a horizontal axis, so that
\begin{align*}
    W_{n-3} &= \begin{pmatrix}
        \bullet & \bullet & & & & \bullet \\
        & \bullet & & & & \\
        & \bullet & & & & \\
        \bullet & & & & & \\
        & & & \bullet & \bullet & \bullet \\
        \bullet & \bullet & \bullet & & &
    \end{pmatrix} &
    W_{n-2} &= \begin{pmatrix}
        & & & \bullet & & \bullet \\
        \bullet & \bullet & \bullet & & & \bullet\\
        & & & & & \bullet\\
        \bullet & & & & \bullet & \\
        \bullet & & & & \bullet & \\
        \bullet & & & & \bullet & \\
    \end{pmatrix} \\
    W_{n-1} &= \begin{pmatrix}
        & & & \bullet & \bullet & \bullet \\
        & & & \bullet & & \\
        & & & \bullet & & \\
        \bullet & \bullet & \bullet & & & \\
        & & & & & \bullet \\
        & & & & \bullet &
    \end{pmatrix} &
    W_n &= \begin{pmatrix}
        & \bullet & & & & \\
        \bullet & & & & & \\
        & & & \bullet & \bullet & \bullet \\
        & & \bullet & & & \\
        & & \bullet & & & \\
        & & \bullet & & &
    \end{pmatrix}.
\end{align*}
Finally, for any $4 < k < n - 3$, the layer $W_k$ is empty.

We first show that introducing a $1$ entry in any of the empty layers of $W$ introduces a copy of $A$. Let $k$ be an empty layer, and suppose we change the entry $o = w_{i,j,k}$ to a $1$. There are four cases, corresponding to the four quarters of layer $k$. We show one in detail; the others are similar.

Suppose that $i, j \geq 4$, so $o$ is in the bottom right quarter of its layer. Then $o$ will correspond to the entry on the second layer of $A$. We define an embedding $\psi_o$ of $A$ into $W$ as follows. $\psi_o(a_{p,q,r}) = w_{p',q',r'}$, where
\begin{align*}
    p' &= \begin{cases} p & p \in \{1,2,3\} \\ i & p = 4 \end{cases} &
    q' &= \begin{cases} q & q \in \{1,2,3\} \\ j & q = 4 \end{cases} &
    r' &= \begin{cases} 4 & r = 1 \\ k & r = 2 \\ n-6+r & r \in \{3,4,5,6\} \end{cases}.
\end{align*}

Now note by inspection that whenever $a_{p,q,r}$ is a $1$ entry in $A$, regardless of the values of $i,j \in \{4,5,6\}$, always $w_{p',q',r'}$ is a $1$ entry in $W$ (after we change $o = \psi_o(a_{4,4,2})$ to a $1$). Thus, when we change $o$ to a $1$, we introduce a copy $\psi_o(A)$ of $A$ in $W$.

Analogous arguments hold for when $o$ lies in the other four quarters of layer $k$. The most notable difference is to which $1$ entry of $A$ the entry $o$ corresponds. When $i,j \leq 3$, i.e. the top left quarter, we have $o = \psi_o(a_{1,1,2})$; when $i \leq 3$ and $j \geq 3$, the top right quarter, $o = \psi_o(a_{1,4,5})$; and when $i \geq 3$ and $j \leq 3$, the bottom left, $o = \psi_o(a_{4,1,4})$.

Now we must check that $W$ does not contain $A$. Suppose that it contains a copy $A^0$, with the embedding $\phi$. The $6$ layers of $A^0$ are some subset of the $8$ nonempty layers of $W$.

Suppose that the layer $A^0_6$ of $A^0$ is $W_n$. It is clear that the first two columns of $A^0_6$ must correspond to the first two columns of $W_n$. Then we search for $\phi(a_{4,2,1})$; this must be $w_{p,2,r}$ for some $p \in \{4,5,6\}$, $r \in \{1,2,3\}$. There are only two possibilities that correspond to $1$ entries: $w_{6,2,1}$ or $w_{5,2,3}$. In the first case, we must have $\phi(a_{3,1,1}) = w_{5,1,1}$, and then $\phi(a_{3,4,6}) = w_{5,q,n}$ for some $q \in \{4,5,6\}$; but none of these $3$ entries is a $1$, ruling out this case. If $\phi(a_{4,2,1}) = w_{5,2,3}$, then we must have $\phi(a_{4,4,2}) = w_{5,q,4}$ for some $q \in \{4,5,6\}$; but again none of these are one entries, so this is impossible. Thus the layer $A^0_6$ is not $W_n$.

Symmetric reasoning (noting that both $A$ and $W$ have twofold rotational symmetry along the horizontal axis) shows that $A^0_1$ is not $W_1$. This leaves only $6$ nonempty layers of $W_n$; thus we in fact see that \[\phi(A^0_k) = \begin{cases} W_{k+1} & k \in \{1,2,3\} \\ W_{k+n-7} & k \in \{4,5,6\} \end{cases}.\] It is fairly clear that, since $\phi(A^0_1) = W_2$, we must have $\phi(a_{1,3,1}) = w_{1,5,2}$ and $\phi(a_{2,4,1}) = w_{2,6,2}$. Also, since $\phi(A^0_6) = W_{n-1}$, we must have $\phi(a_{3,4,6}) = w_{5,6,(n-1)}$, and $\phi(a_{4,3,6}) = w_{6,5,(n-1)}$. But then $\phi(a_{3,1,1}) = w_{5,q,2}$, for some $q \in \{1,2,3\}$; but none of these are $1$ entries, so this is impossible. Thus $A$ is not contained in $W$.

Now, we may define $\Tilde{W}$ as follows: starting from pattern $W$, change zeroes to ones until $\Tilde{W}$ is saturating for $A$. Since changing any of the zeroes in the empty layers would introduce a copy of $A$, we know that $\Tilde{W}$ still has the $n - 8$ empty layers in the middle that $W$ had. Then since $\Tilde{W}$ is a $6 \times 6 \times n$ matrix with at most $8$ nonempty layers, its weight is at most $6 \cdot 6 \cdot 8 = 288$. Furthermore, $\Tilde{W}$ is $A$-saturating. Together this shows $\sat(6,6,n; A) \leq 288$, so $\sat(6,6,n; A) = O(1)$ as desired.
\end{proof}

This is the first known example of a pattern of dimension greater than $2$ with any saturation function bounded. Our proof is analogous to the use of a witness in the $2$-dimensional case; we proceed by presenting what is essentially a 3-dimensional ``(vertical) witness''; we first show that it has an ``expandable layer'', and then that it does not contain the forbidden pattern.

We also extend a theorem of Geneson and Tsai \cite{gt23} regarding the multidimensional semisaturation function. They showed that, for any $d$-dimensional pattern $P$, the function $\ssat(n,P) = \Theta(n^k)$ for some $k \in \{0,\dots,n - 1\}$. Generalizing both our Proposition \ref{prop:vhssat} and this theorem of \cite{gt23}, we show that the same holds when we fix certain dimensions and let the others go to infinity---that is, we have

\begin{theorem}
    For any $d$-dimensional $p_1 \times \dots p_d$ pattern $P$, any integer $\ell \in [d]$, and any fixed integers $m_1, \dots, m_\ell$, the semisaturation function \[\ssat(m_1,\dots,m_\ell,n,\dots,n; P) = \Theta(n^k)\] for some $k \in [d - \ell]$.

    In particular, $k$ is the least positive integer satisfying    \begin{itemize}[label=$(\ast)$]
        \item If $f$ is a face of $P$ with dimension at least $k$, then $f$ has a one entry $o$ such that, if $g$ is a cross-section of $P$ with $|(C_g \setminus C_f) \setminus [\ell]| \geq k + 1$ containing $o$, then $o$ is the only $1$ entry in $g$.
    \end{itemize}
\end{theorem}

\begin{proof}
    Fix a $p_1 \times \dots p_d$ pattern $P$, $\ell \in [d]$, and $m_1, \dots, m_\ell$. Let $k$ be the least positive integer satisfying $(\ast)$ above. Take any $n > \max\{p_i : i \in [d]\}$.
    
    For brevity write $m_i = n$ for integers $i \in [\ell + 1, n]$), so that the expression $m_1 \times \cdots \times m_\ell \times n \times \cdots \times n$ may be written as $m_1 \times \cdots \times m_d$.
    
    We will show that $\ssat(m_1,\dots,m_\ell,n,\dots,n; P) = \Theta(n^k)$. First \textbf{(Step 1)}, we construct \textbf{(Step 1a)} an $m_1 \times \cdots \times m_d$ pattern $M$; we show \textbf{(Step 1b)} $M$ is semisaturating, and then \textbf{(Step 1c)} we see that $M$ has $O(n^k)$ one entries. Then \textbf{(Step 2)}, we show that if $M$ is a $P$-semisaturating matrix, then $w(M) = \Omega(n^k)$, by double counting the number of zero entries of $M$. First \textbf{(Step 2a)}, using the $P$-semisaturating property and $(\ast)$, we give a precise connection between zero entries and one entries in $M$; we then find \textbf{(Step 2b)} that $M$ can have only a limited number of zero entries per one entry. Then \textbf{(Step 2c)}, we count the total number of entries. Combining these \textbf{(Step 2d)}, we get a lower bound of $\Omega(n^k)$ for the number of one entries $w(M)$.
    
    \textbf{Step 1.} We first show that $\ssat(m_1,\dots,m_\ell,n,\dots,n; P) = O(n^k)$.
    
    \textbf{Step 1a.} Construct the $m_1 \times \cdots \times m_d$ matrix $M$ as follows. For the entry $o = (o_1, \dots, o_d)$, define the set
    \[F_o = \{i \in [\ell + 1, d] : o_i \in [p_i, m_i - p_i + 1]\}.\]
    Then make $o$ a $1$ entry in $M$ if and only if $|F_o| \leq k$. We will show $M$ is semisaturating for $P$.

    \textbf{Step 1b.} Fix an entry $o = (o_1, \dots, o_d)$ that is zero in $M$, and flip $o$ to a one entry to make the new matrix $M_o$. Since $o$ was originally a zero entry, we must have $|F_o| > k$. Consider $i \in [d] \setminus F_o$; we have either $o_i \leq p_i - 1$ or $o_i \geq m_i - p_i + 2$. Define \[q_i = \begin{cases}
        1 & o_i \leq p_i - 1 \\ p_i & o_i \geq m_i - p_i + 2
    \end{cases}.\] Consider the face $f$ of $P$ where $C_f = [d] \setminus F_o$, and the $i$-th coordinate is fixed to $q_i$. There is a one entry $o' = (o'_1, \dots, o'_d)$ in $f$ satisfying the conditions of $(\ast)$, by hypothesis on $k$.

    We build the submatrix $P_o$ of $M_o$ as follows. For any $i \in [d]$, we let the $i$-th coordinate range over $\{1,\dots,o'_i - 1\} \cup \{o_i\} \cup \{m_i - p_i + o'_i + 1, \dots, m_i\}$; then $P_o$ is the $p_1 \times \cdots \times p_d$ matrix built by taking these possible coordinates in each dimension $i$. We claim that $P_o$ corresponds to a copy of $P$ in $M_o$, that is, that whenever $u'$ is a one entry in $P$, the corresponding entry $u$ in $P_o$ is one.

    Suppose that $u'$ is a one entry in $P$, corresponding to the entry $u = (u_1, \dots, u_d)$ in $P_o$. Notice that, because of the restrictions on possible coordinates in $P_o$, $u_i \in [p_i,m_i - p_i + 1]$ only if $u_i = o_i$. Then \[F_u = \{i \in F_o: u_i = o_i\}.\]
    
    If $|F_u| \leq k$, then certainly $u$ is a one entry by construction of $M$. Otherwise, $|F_u| \geq k + 1$.
    
    Take the cross section $g$ of $P_o$ where $C_g = F_u$ and each coordinate is fixed to $o_i$; then both $o$ and $u$ lie in $g$. Then, the corresponding cross section $g'$ in $P$ also has $C_{g'} = C_g= F_u$, and this time each coordinate is fixed to $o'_i$. Since both $o$ and $u$ were in $g$, also both $o'$ and $u'$ are in $g'$. Furthermore, $C_g \subseteq F_o$, so $C_g \cap C_f = \emptyset = C_g \cap [\ell]$. Then $|(C_g \setminus C_f) \setminus [\ell]| = |C_g| \geq k + 1$. Then, by $(\ast)$, $o'$ is alone in $g$. But we assumed that $u'$ was a one entry as well; thus $o' = u'$. Then the corresponding entry $o = u$, so $u$ is also a one entry.

    Thus, every one entry of $P$ corresponds to a one entry in $P_o$, so $P_o$ is a copy of $P$ in $M_o$. Thus $M$ is semisaturated for $P$ as desired.

    \textbf{Step 1c.} Now, we count $w(M)$. First, let $S \subseteq [\ell + 1, d]$ with $|S| = k$, and define the submatrix of $M$
    \[E_S = \{o : F_o \subseteq S\} = \{o : (\forall i \in [d] \setminus S)\: o_i \notin [p_i, m_i - p_i + 1]\}.\]
    Observe that each one entry of $M$ is in $E_S$ for some $S$, and there are $\binom{d - \ell}{k}$ possibilities for $S$. Then we have
    \[w(M) \leq \binom{d - \ell}{k} |E_S|\]
    where the notation $|E_S|$ is the total number of entries in the matrix $E_S$. We can count $|E_S|$ based on the possible coordinates of each entry. Let $o = (o_1, \dots, o_d) \in E_S$. Then, if $i \notin S$ and $i > \ell$, we have $i \notin [p_i, m_i - p_i + 1]$. Then we have the following for every $i$:
    \[o_i \in \begin{cases}
        [m_i] & i \leq \ell \\
        [n] & i \in S \\
        [1, p_i - 1] \cup [m_i - p_i + 2, m_i] & i \notin S \land i > \ell
    \end{cases}\]
    Then, multiplying these possible coordinates,
    \begin{align*}
        |E_S| &= \left(\prod_{i = 1}^\ell m_i \right)\left(\prod_{i \in S} n\right)\left(\prod_{i \notin S, i > \ell} 2(p_i - 1)\right) \\
        &= \left(\prod_{i = 1}^\ell m_i \right)\left(\prod_{i \notin S, i > \ell} 2(p_i - 1)\right) n^k \\
        &\leq \left(\prod_{i = 1}^\ell m_i \right) (2p - 2)^{d - \ell - k} n^k
    \end{align*}
    where $p = \max \{p_i : i \in [\ell + 1, d]\}$. Notice that this bound does not depend on $S$.  Thus,
    \[w(M) \leq \binom{d - \ell}{k} |E_S| \leq \binom{d - \ell}{k} \left(\prod_{i = 1}^\ell m_i \right) (2p - 2)^{d - \ell - k} n^k = O(n^k).\]
    Thus $\ssat(m_1,\dots,m_\ell,n,\dots,n; P) = O(n^k)$ as desired.

    \textbf{Step 2.} Now we must show the lower bound: $\ssat(m_1,\dots,m_\ell,n,\dots,n; P) = \Omega(n^k)$. Let $M$ be any $m_1 \times \cdots \times m_d$ matrix that is semisaturating for $P$.
    
    \textbf{Step 2a.} By maximality of $k$, there is some face $f$ of $P$ of dimension $k - 1$ violating $(\ast)$. Call a zero entry and a one entry of $M$ \textit{connected} if they lie in a shared cross-section $g$ of $M$ such that $|(C_g \setminus C_f) \setminus [\ell]| = k$.
    
    Let $z$ be a zero entry of $M$. If $z$ is changed to a $1$, making the matrix $M_z$, then by semisaturation there is a copy $P_z$ of $P$ in $M_z$ so that $z$ is one of the one entries of $P_z$. The violation of $(\ast)$ is exactly the statement that there is some choice of $C_g$ so that $z$ is not the only one entry in the cross-section $g$ and $|(C_g \setminus C_f) \setminus [\ell]| \geq k$; that is, $z$ is connected in $P_z$ to some one entry $o$. That is, any zero entry of $M$ must be connected to a one entry in $M$.
    
    \textbf{Step 2b.} Now fix $o = (o_1,\dots,o_d)$ and suppose that $z = (z_1,\dots,z_d)$ is connected to $o$; they lie in a shared cross-section $g$ per above. We can assume $g$ is large enough that $C_g \cap C_f = C_g \cap [\ell] = \emptyset$ (otherwise, pick a larger cross-section containing $g$). Then $|C_g| = k$, and, since $|([d] \setminus C_f) \setminus [\ell]| = d - k - \ell + 1$, there are $\binom{d - k - \ell + 1}{k}$ possibilities for $C_g$, and hence $\binom{d - k - \ell + 1}{k}$ possibilities for $g$. Each possible $g$ is $(d - k)$-dimensional and hence contains $m_1m_2\cdots m_\ell n^{d - k - \ell}$ entries, so $o$ is connected to at most $\binom{d - k - \ell + 1}{k} m_1m_2\cdots m_\ell n^{d - k - \ell}$ entries.

    Then, since each zero entry is connected to some one entry, the number of zero entries is at most
    \[\binom{d - k - \ell + 1}{k}\left( \prod_{i=1}^\ell m_i \right) n^{d - k - \ell} w(M).\]
    
    \textbf{Step 2c.} But the number of zero entries is also given more simply by counting the total entries. Since there are $\left(\prod_{i=1}^\ell m_i\right)n^{d - \ell}$ total entries in the matrix, the number of zero entries is exactly $\left(\prod_{i=1}^\ell m_i\right)n^{d - \ell} - w(M)$. 
    
    \textbf{Step 2d.} Using the observation of \textbf{Step 2c}, the bound from \textbf{Step 2b} can be rearranged to give:
    \begin{align*}
        \left(\prod_{i=1}^\ell m_i\right)n^{d - \ell} - w(M) &\leq \binom{d - k - \ell + 1}{k} \left( \prod_{i=1}^\ell m_i \right) n^{d - k - \ell} w(M) \\
        w(M) &\geq \frac{\left(\prod_{i=1}^\ell m_i\right)n^{d - \ell}}{1 + \binom{d - k - \ell + 1}{k}\left(\prod_{i=1}^\ell m_i\right)n^{d - k - \ell}} \\
        &= \Omega(n^k).
    \end{align*}
    Then, any $m_1 \times \cdots \times m_d$ matrix $M$ satisfies $w(M) = \Omega(n^k)$. Thus, \[\ssat(m_1,\dots,m_\ell,n,\dots,n; P) = \Omega(n^k),\] as desired.
\end{proof}

In the remainder of this section, we generalize the column-adding and row-adding operations for $\sat(n, P)$ from Section~\ref{sec:insert} into layer-adding operations for $\sat(n, P, d)$. Then, we prove a sharp upper bound on $\sat(n, P, d)$ with respect to the dimensions of $P$.

\begin{thm}\label{thm_dcol_add_op} Let $P$ be any $d$-dimensional pattern, and let $P'$ be any pattern derived from $P$ by adding a single layer before the first layer of $P$. Then \[\sat(n_1,n_2,\dots,n_d,P') \le n_1 n_2 \dots n_{d-1} + \sat(n_1,n_2,\dots,n_{d-1},n_d-1,P).\] \end{thm} 

\begin{proof} Let $M$ be an $n_1 \times n_2 \times \dots \times n_{d-1} \times (n_d-1)$ matrix that is saturating for $P$. Consider the matrix $M'$ derived from $M$ by adding a layer of all $1$ entries before $M$. We claim that $M'$ is saturating for $P'$. Suppose for contradiction that $M'$ contains $P'$. Removing the leftmost layer of $P'$, we obtain a copy of $P$ in $M'$; furthermore this copy of $P$ does not include the first layer of $M'$, so $M$ also contains a copy of $P$, which is a contradiction. To see that $M'$ is saturating, we now show that adding any $1$ to $M'$ introduces a copy of $P'$. Since the first layer of $M'$ is all $1$, any $1$ entry must be added outside the first layer. Then this $1$ entry is in $M$, so it is part of a copy of $P$ in $M$; this copy of $P$ together with the first layer of $M'$ forms a copy of $P'$ in $M'$ as desired. Thus, \[\sat(n_1,n_2,\dots,n_d,P') \le n_1 n_2 \dots n_{d-1} + \sat(n_1,n_2,\dots,n_{d-1},n_d-1,P).\] \end{proof} 

 Note that the same bound also holds for adding a single layer after the last layer of $P$. Moreover, the same argument also applies to adding $i$-layers for any $i < d$. These variants also hold for the next corollary.

\begin{cor}\label{cor_empty_dcol_add_op} Let $P$ be any pattern, and let the pattern $P'$ be derived from $P$ by adding an empty
layer before the first layer of $P$. Then, \[\sat(n_1,n_2,\dots,n_d,P') = n_1 n_2 \dots n_{d-1} + \sat(n_1,n_2,\dots,n_{d-1},n_d-1,P).\] \end{cor} 

\begin{proof} By Theorem~\ref{thm_dcol_add_op}, we have
\[\sat(n_1,n_2,\dots,n_d,P') \le n_1 n_2 \dots n_{d-1} + \sat(n_1,n_2,\dots,n_{d-1},n_d-1,P).\] For the other direction, now let $M$ be any
$n_1 \times n_2 \times \dots \times n_d$ matrix that is saturating for $P'$. If the first layer of $M$ has any $0$ entries, then changing them to a $1$ cannot introduce a copy of $P'$, since the first layer of $P'$ is empty. Thus the first layer of $M$ must have all ones. We claim that the matrix $M'$ derived from $M$ by deleting the first layer is saturating for $P$. Clearly if $M'$ contains a copy of $P$, then this copy of $P$ with the first layer of $M$ added forms a copy of $P'$ in $M$, so $M'$ cannot contain a copy of $P$. If a $0$ entry in $M'$ is changed to a $1$, we know that this creates a copy of $P'$ in $M$; then this also introduces a copy of $P$ in $M$ which cannot include the first layer, which is a copy of $P$ in $M'$. Thus, $M'$ is saturating for $P$, so \[\sat(n_1,n_2,\dots,n_{d-1},n_d-1,P') \le \sat(n_1,n_2,\dots,n_d,P)-n_1 n_2 \dots n_{d-1}.\] Combining the two directions of inequality, the result follows. \end{proof}

In the next result, we show that $\sat(n, P, d)$ is maximized over all $d$-dimensional 0-1 matrices $P$ of dimensions $k_1 \times k_2 \times \dots \times k_d$ when $P$ has a single entry equal to one and zeroes elsewhere. Tsai \cite{tsai23} proved that \[\sat(n,P,d) \le n^d - \prod_{i = 1}^d (n - k_i + 1)\] for all $d$-dimensional 0-1 matrices $P$ of dimensions $k_1 \times k_2 \times \dots \times k_d$. The construction in \cite{tsai23} implies a more general bound of \[\sat(n_1,n_2,\dots,n_d,P,d) \le \prod_{i = 1}^d n_i - \prod_{i = 1}^d (n_i - k_i + 1)\] for $n_i \ge k_i-1$. In the next result, we show that this bound is sharp.

\begin{thm}
    For $n_i \ge k_i-1$, the maximum possible value of the saturation function $\sat(n_1,n_2,\dots,n_d,P,d)$ over all $d$-dimensional 0-1 matrices $P$ of dimensions $k_1 \times k_2 \times \dots \times k_d$ is \[\prod_{i = 1}^d n_i - \prod_{i = 1}^d (n_i - k_i + 1).\]
\end{thm}

\begin{proof}
    The upper bound follows from the construction in \cite{tsai23}. For the lower bound, consider the 0-1 matrix $P$ of dimensions $k_1 \times k_2 \times \dots \times k_d$ with a single one in the entry with all coordinates maximal, and zeroes elsewhere. Clearly for $n_i \ge k_i-1$, any $P$-saturated $d$-dimensional 0-1 matrix with dimensions of length $n_1 \times n_2 \times \dots \times n_d$ must have all ones in its first $k_i-1$ $i$-layers for each $1 \le i \le d$. This implies the lower bound.
\end{proof}

\section{Witness Graphs}\label{sec:witness_graphs}

In this section, we introduce a combinatorial construction which captures much of the data of a witness. We make a few general observations and then use this tool to resolve two specific $2 \times 4$ patterns.

\begin{definition}
    Let $P$ be a pattern and $W$ any horizontal (vertical) witness for $P$. A \textit{witness graph} $G^j(W)$ (resp. $G_j(W)$) associated to $W$ for the expandable column (resp. row) $j$ is a directed graph whose vertices are the rows (resp. columns) of $W$, with an edge from row $a$ to row $b$ if $a \neq b$ and row $b$ is contained within the submatrix $P_a$. When the expandable column $j$ is understood, we write simply $G(W)$.
\end{definition}

Note that this is not defined uniquely for a given $W$; there may be multiple witness graphs. In particular, there is a witness graph corresponding to each possible choice of copy $P_i$ of $P$ in $W$.

We will state results for the witness graph of horizontal witnesses. These all similarly apply to vertical witnesses.

\begin{definition}
    Let $P$ be any pattern. A witness $W$ for $P$ is minimal if no submatrix of $W$ is a witness for $P$.
\end{definition}

This definition is useful because of the following observation.

\begin{obs}
    Any witness $W$ for $P$ can be refined into a minimal witness.
\end{obs}

This means in particular that, in order to show that no witness exists, it suffices to show that no minimal witnesses exist; we may restrict our search to the case of minimal witnesses.

\begin{prop} \label{prop:gw_corr}
    If $W'$ is a submatrix of $W$ that is also a witness for $P$ (with the same expandable column), then any witness graph $G(W')$ is a nonempty subgraph of a witness graph $G(W)$ such that if row $a \in G(W')$, then every row $b \in G(W)$ with an edge connecting $a$ to $b$ is also in $G(W')$. Conversely, for any such subgraph of a witness graph $G(W)$, the corresponding submatrix is a witness for $P$.
\end{prop}

\begin{figure}
    \begin{align*}
        P &= \begin{pmatrix}
            & \bullet & & \\
            & & & \bullet \\
            \bullet & & & \\
            & & \bullet &
        \end{pmatrix} &
        W &= \begin{pmatrix}
            & & & \bullet & & \\
            & \bullet & & & & \bullet \\
            \cdot & \cdot & \cdot & \cdot & \cdot & \cdot \\
            \bullet & & & & \bullet & \\
            & & \bullet & & &
        \end{pmatrix}
    \end{align*}\\
    \centering
    \begin{tikzpicture}[scale=0.5, -stealth]
        \node (1) {1};
        \node (2) [right of=1] {2};
        \node (3) [right of=2] {3};
        \node (4) [right of=3] {4};
        \node (5) [right of=4] {5};
        \node (6) [right of=5] {6};
    
        \path[every node/.style={font=\sffamily\small}]
        (1) edge[bend left] node [right] {} (4)
        (1) edge[bend left] node [right] {} (5)
        (1) edge[bend left] node [right] {} (6)
        (2) edge[bend left] node [right] {} (4)
        (2) edge[bend left] node [right] {} (5)
        (2) edge[bend left] node [right] {} (6)
        (3) edge[bend left] node [right] {} (4)
        (3) edge[bend left] node [right] {} (5)
        (3) edge[bend left] node [right] {} (6)
        (4) edge[bend left] node [left] {} (1)
        (4) edge[bend left] node [left] {} (2)
        (4) edge[bend left] node [left] {} (3)
        (5) edge[bend left] node [left] {} (1)
        (5) edge[bend left] node [left] {} (2)
        (5) edge[bend left] node [left] {} (3)
        (6) edge[bend left] node [left] {} (1)
        (6) edge[bend left] node [left] {} (2)
        (6) edge[bend left] node [left] {} (3);
    \end{tikzpicture}
    \caption{An example of $G_3(W)$ for pattern $P$ and witness $W$. This is isomorphic to the complete bipartite directed graph $K_{3,3}$.}
    \label{fig:q1_witness_graph}
\end{figure}
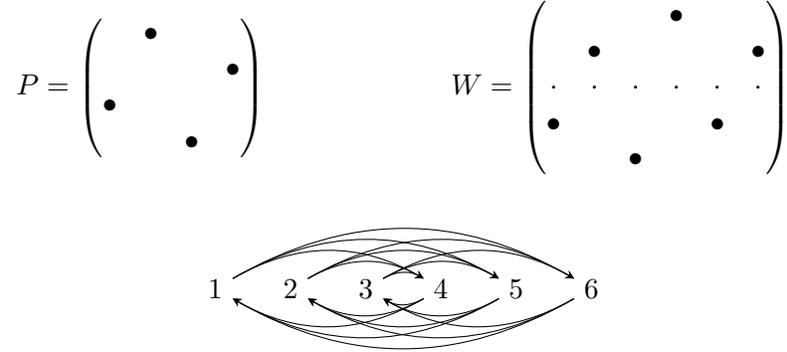

\begin{proof}
    Let $W$ be a witness for $P$. For the forward direction, let $W'$ be a submatrix of $W$ that is also a witness for $P$. Pick patterns $P_i$ for each row $i$ in $W'$ and extend this to a choice of $P_i$ for all rows in $W$; in this way we can ensure that for each row $i$ of $W'$, $P_i$ lies entirely within $W'$. On the graph $G(W)$ given by this choice of $P_i$, the rows of $W'$ correspond to a subgraph. To see the claim, let $a$ be any row of $W'$; if there is an edge from row $a$ to row $b$, then $b$ is in $P_i$, so by our choice of $P_i$ above we know $b$ is in $W'$, so $b$ is in $G(W')$.

    For the other direction, let $S$ be a subgraph of $G(W)$ satisfying the conditions given. Take the submatrix $W'$ of $W$ which contains the rows in $S$. Certainly $W'$ avoids $P$. To see that $W'$ is a witness, note that the condition on $S$ guarantees that, for each row $i$ of $W'$, the entire matrix $P_i$ lies in $W'$, so $W'$ will still have an expandable column.
\end{proof}

The converse of this gives the following useful corollary.

\begin{cor}
    If $W$ is a minimal witness for $P$, then $G(W)$ is connected, and every vertex has at least one edge leading into it.
\end{cor}

\begin{proof}
    If $G(W)$ is not connected, then take a subgraph corresponding to one of its connected components; this subgraph satisfies the conditions of Proposition \ref{prop:gw_corr}, and so it corresponds to a proper submatrix of $W$ which is a witness. Then $W$ is not minimal.

    Similarly, if $G(W)$ has a vertex without an edge leading into it, we may consider the subgraph obtained by deleting this vertex to show $W$ is not minimal.
\end{proof}

We also have the following simple facts.

\begin{prop}
    If $P$ has $k$ rows and $W$ is a witness for $P$, then each vertex of $G(W)$ has $k - 1$ edges leading out of it.
\end{prop}

\begin{proof}
    Each edge corresponds to a row of $P$, except for one which would be a self-loop (which the definition excludes). There are $k$ rows of $P$, so $k - 1$ edges leading out of any given row.
\end{proof}

\begin{prop}
    For any witness $W$, any witness graph $G(W)$ is cyclic; in particular, there exists a path from any point to a point on a cycle.
\end{prop}

\begin{proof}
    Build a trail starting at any vertex in $P$. Since all vertices have at least one edge leading out, this trail can always be extended. By the pigeonhole principle, it must eventually repeat a vertex, which gives a cycle.
\end{proof}

We note that known examples  of witness graphs are often bipartite. The example in Figure \ref{fig:q1_witness_graph} is a complete bipartite directed graph. The construction in \cite{berendsohn21}, which generalizes this example to show Theorem \ref{thm:q1like}, will always give a witness with witness graph equal to a complete bipartite directed graph. Furthermore, the witness construction in \cite{berendsohn23} for permutation matrices also always gives a bipartite graph, although in this case the graph is not complete. Furthermore, many of the particular examples treated in Section \ref{sec:patterns} have witnesses with bipartite witness graphs, including those covered by Theorem~\ref{thm:q2like}. However, not all witness graphs are bipartite, as the following example illustrates.

\begin{example}
    The pattern \[P = \begin{pmatrix}
        \bullet & & \bullet & \bullet \\
        \bullet & \bullet & & \bullet
    \end{pmatrix}\] has a witness \[W = \begin{pmatrix}
        & \bullet & \bullet & \cdot & \bullet & \bullet \\
        & \bullet & \bullet & \cdot & \bullet & \\
        & & \bullet & \cdot & \bullet & \bullet \\
        \bullet & \bullet & \bullet & \cdot & \bullet & \\
        \bullet & & \bullet & \cdot & \bullet & \bullet
    \end{pmatrix}\] which has $G(W)$ equal to a cycle on five vertices (which is not a bipartite graph).
\end{example}

In the remainder of this section, we apply some of the facts above to show that the particular patterns
\[S = \begin{pmatrix}
    & \bullet & & \bullet \\
    \bullet & & \bullet &
\end{pmatrix}, \hspace{0.5 cm} 
S' \begin{pmatrix}
    & \bullet & & \bullet \\
    \bullet & & \bullet & \bullet
\end{pmatrix}\] have no horizontal witnesses. This proof includes some similar ideas to that of Proposition \ref{prop:q3_like}, but the latter is not phrased in the language of witness graphs or minimal witnesses.

\begin{lem}\label{lem:witness_graph_cycle}
    Suppose $P$ is any $2 \times \ell$ pattern, and $W$ is a minimal witness for $P$. Then any witness graph $G(W)$ is a cycle graph.
\end{lem}

\begin{proof}
    Each vertex of $G(W)$ has exactly one edge out of it, and at least one edge into it. There are exactly as many edges as vertices, so each vertex must in fact have exactly one edge into it. Then there exists a unique path beginning at each vertex; this path must eventually cycle. Since only one edge can lead into each vertex, the first vertex which cycles must in fact lie at the beginning of the path, so each vertex lies on a cycle. Then because $G(W)$ is connected, $G(W)$ is in fact one cycle.
\end{proof}

For a pattern $P$ as in \ref{lem:witness_graph_cycle}, define the map $\psi_P$ which takes each row $i$ to the unique row to which it is connected by an edge. Iterations of $\psi_P$ go along the cycle $G(W)$. Note that in fact this lemma, and this map, are applicable in general to any pattern with two rows.

\begin{lem}\label{lem:p_all_1s}
    With $S$, $S'$ as above, suppose that $W$ is a minimal witness for $S$ or $S'$ with expandable column $j$. Then columns $j - 1$ and $j + 1$ are all-ones.
\end{lem}

\begin{proof}
    We only show the case for $S$; the proof applies verbatim to $S'$.
    
    Clearly $j$ is not the first column of $W$ (since the top left corner must be a $1$ entry in any $S$-saturated matrix). Since $W$ is minimal, column $j - 1$ is nonempty; then some entry of column $j - 1$ is a $1$ entry. Let this entry be in the $i$th row.
    
    Suppose first that $i$ is in the top row of $S_i$. Then if the first column of $P_i$ is to the left of $j - 1$, we replace column $j$ of $P_i$ with $j - 1$ and obtain a copy of $S$ in $W$ itself, which is a contradiction. Then the first column of $S_i$ must be $j - 1$, so in particular the entry in row $\psi(i)$ must be a one entry.

    Similarly, if $i$ is in the bottom row of $S_i$, we also see that the entry in row $\psi(i)$ is a one entry; here we use the second column of $P_i$.

    So in either case, the entry in row $\psi(i)$ is also $1$. Iterating $\psi$ covers every row of $W$, by Lemma~\ref{lem:witness_graph_cycle}; thus, applying this argument inductively, every entry of column $j - 1$ must be $1$.

    A completely symmetric argument shows the claim for column $j + 1$ as well.
\end{proof}

\begin{figure}
    \begin{align*}
    &\begin{pmatrix}
        & \circ & \cdot & \bullet & & \circ \\
        & \bullet & \cdot & \bullet & & \\
        & & \vdots & & & \\
        & \bullet & \cdot & \bullet & & \\
        \circ & \bullet & \cdot & \circ & & \\
    \end{pmatrix}
    & & \begin{pmatrix}
        & \circ & \cdot & \bullet & & \circ \\
        & & \vdots & & & \\
        \circ & \bullet & \cdot & \circ & & \circ\\
        & & \vdots & & & \\
        & \bullet & \cdot & \bullet & & \\
    \end{pmatrix}
    \end{align*}
    \caption{Finding patterns $S$ and $S'$ in candidate witnesses.}
    \label{fig:2x4s}
\end{figure}

\begin{thm}
    For the patterns $S$ and $S'$ as above and any $m_0$, we have $\sat(m_0,n,S) = \Theta(n)$ and $\sat(m_0,n,S') = \Theta(n)$.
    \label{thm:2x4s}
\end{thm}

\begin{proof}
    First, we show the case for $S$. By Lemma \ref{lem:p_all_1s}, columns $p - 1$ and $p + 1$ are all-ones. Let $a$ be the last column of $S_1$; $a$ must be at least $p + 2$, and the entry $w_{1a} = 1$. Similarly let $b$ be the first column of $S_m$; $b \leq p -2$, and the entry $w_{mb} = 1$ as well. Then the submatrix with columns $\{b, p-1, p+1, a\}$ and rows $\{1, m\}$ is a copy of $S$ in $W$.

    For $S'$, our argument is broadly similar. Again, columns $p - 1$ and $p + 1$ are all-ones. Let $a$ be the last column of $S'_1$, and $r = \psi_{S'}(1)$. The column $a$ must be at least $p + 2$, and both entries $w_{1a} = w_{ra} = 1$. Let $b$ be the first column of $S'_r$; then $b \leq p - 2$, and the entry $w_{rb} = 1$ as well. Then form $S'_0$ by taking columns $\{b, p-1, p+1, a\}$ and rows $\{1,r\}$; this matrix is a copy of $S'$.
\end{proof}

\section{Linear Program}\label{sec:lin_prog}
Now, we present a computational approach that models saturation using a binary integer linear program.
\begin{thm}
    Let $P$ be an $m\times n$ pattern with weight $w(P)$. Let $\mathcal P$ be the set of all sets of indices which form a copy of $P$ in the $m\times n$ all-ones matrix $\mathbf1_{m\times n}$ and $A_{i,j}=\{X\in \mathcal P:(i,j)\in X\}$ be all copies of $P$ containing $(i,j)$. Notice that $\lvert P'\rvert=w(P)$ for all $P'\in\mathcal P$. Define variables $x_{i,j}$ for $1\le i\le m$ and $1\le j\le n$ and $y_{P'}$ for each $P'\in \mathcal P$. For notational convenience, set $f(P')=\sum_{(i,j)\in P'}x_{i,j}$. Suppose $x_{i,j}$ and $y_{P'}$ satisfy
    \begin{alignat}3
        x_{i,j}&\in\{0,1\}&\forall(i,j)\in[m]\times[n]\\
        y_{P'}&\in\{0,1\}&\forall P'\in\mathcal P\\
        \label{eqn:not_contain}
        f(P')&<w(P)&\forall P'\in\mathcal P\\
        \label{eqn:yp'_cond_1}
        (w(P)-1)y_{P'}-f(P')&\le0&\forall P'\in\mathcal P\\
        \label{eqn:yp'_cond_2}
        y_{P'}-f(P')&\ge-w(P)+2&\forall P'\in\mathcal P\\
        \label{eqn:changeone}
        \sum_{P'\in A_{i,j}}y_{P'}&>0&\forall(i,j)\in[m]\times[n]\\
        \label{eqn:lastcheck}
        x_{i,j}&=1&\forall(i,j)\in[m]\times[n],A_{i,j}=\varnothing.
    \end{alignat}
    Then $\min\sum_{i,j}x_{i,j}=\sat(m,n,P)$.
\end{thm}
\begin{proof}
    We claim that if we have a matrix $M$ with $m_{ij}=x_{i,j}$, then $M$ is saturating for $P$ and $y_{P'}=1$ if and only if $m_{i,j}=1$ for exactly $w(P)-1$ pairs $(i,j)\in P'$ for all $P'\in\mathcal P$. This suffices because $\sat(m,n,P)$ is defined to be the minimum possible weight of such a saturating matrix $M$.

    Any saturating matrix $M$ must not contain $P$. This is satisfied by equation \eqref{eqn:not_contain} because it states that for every possible copy $P'$ of $P$, less than $w(P)$ of those entries are actually $1$ in $M$.

    Fix some $P'$. The condition that $y_{P'}=1$ if and only if all but one entry of $P'$ is in $M$ is satisfied through equations \eqref{eqn:yp'_cond_1} and \eqref{eqn:yp'_cond_2}. Let $f(P')=\sum_{(i,j)\in P'}x_{i,j}$, the number of entries of $P'$ present in $M$. If $y_{P'}=1$, then \eqref{eqn:yp'_cond_1} implies that $f(P')\ge w(P)-1$ while \eqref{eqn:yp'_cond_2} implies $w(P)-1\ge f(P')$, so $f(P')=w(P)-1$ as claimed. Now if $y_P'=0$, then \eqref{eqn:yp'_cond_1} and \eqref{eqn:yp'_cond_2} become $f(P')\ge0$ and $f(P')\le w(P)-2$, which implies that less than $w(P)-1$ entries of $P'$ are present in $M$. This proves the condition on $y_{P'}$.

    Finally, we show that changing any $0$ to a $1$ in $M$ introduces a copy of $P$. If some entry of $M$ isn't part of any $P'\in\mathcal P$, then it must be $1$ or else changing that entry to $1$ wouldn't introduce a copy of $P$. This is encoded in \eqref{eqn:lastcheck}. Finally, \eqref{eqn:changeone} says that for every $(i,j)$, at least one of the $P'$ that contain $(i,j)$ has $w(P)-1$ entries in $M$. This means that changing the $(i,j)$ entry to $1$ will introduce a copy of $P$ because then some $P'$ would have $w(P)$ entries in $M$, meaning $M$ contains $P$. 
\end{proof}
An implementation of this model using Sage can be found in \cite{code}. This model is useful for determining the exact values of saturation functions, which is an unsolved problem for most forbidden patterns. As an example of using the model to evaluate exact values of the saturation function, we turn our attention to a particular pattern,
\[P=\begin{pmatrix}
    \bullet & \bullet &  \\
     & \bullet & \bullet \\
     &  & \bullet
\end{pmatrix}\]
Using the implementation in \cite{code}, we find the following values:
\begin{table}[h!]
\centering
\begin{tabular}{c|c}
    $m,n$&$\sat(m,n,P)$\\\hline
    $3,3$&8\\
    $3,4$&10\\
    $4,4$&12\\
    $4,5$&14
\end{tabular}
\caption{Values of the saturation function for small $(m,n)$.}
\end{table}

From the data, we conjecture that $\sat(m,n,P)=2m+2n-4$ for $m,n\ge3$. In fact, we obtain the same values from any pattern of the form
\[T=\begin{pmatrix}
    \bullet&*&*\\
    &\bullet&*\\
    &&\bullet
\end{pmatrix}.\]
If we were to generalize to larger sizes of $k\times k$ upper triangular matrices with a full diagonal, we would hypothesize that $\sat(m,n,P)=(k-1)(m+n-k+1)$. Note that this is the general upper bound on the $\sat$ function given by Fulek and Keszegh \cite{FK}. However, the following $4\times 4$ pattern is a counterexample.
\[
    P=\begin{pmatrix}
        \bullet&&&\bullet\\
        &\bullet&&\\
        &&\bullet&\\
        &&&\bullet
    \end{pmatrix}
\]
In this case, the linear program yields $\sat(5,5,P)=20$. Therefore, we ask when is this upper bound sharp, besides in the known case of the identity matrix from Proposition \ref{thm:idmat}?
\begin{problem}
    Let $T$ be any $k\times k$ triangular 0-1 matrix with all ones on the diagonal. If $k=3$, we conjecture that $\sat(m,n,T)=2m+2n-4$ for $m,n\ge3$. For $k>3$, what conditions need to be placed to ensure that $\sat(m,n,T)=(k-1)(m+n-k+1)$ for $m,n\ge k$?
\end{problem}

\section{Specific Patterns}
\label{sec:patterns}

We conclude by resolving the saturation of various small matrices and several infinite families of matrices. Our results are sufficient to completely classify $4 \times 4$ matrices in terms of their saturation function, in addition to many others.

First, we observe that existing results suffice to resolve many small patterns. Recall the following family of matrices from \cite{berendsohn21}.

\begin{thm}[\cite{berendsohn21}]
    Let $P$ be any pattern satisfying all of the following.
    \begin{enumerate}
        \item $P$ contains exactly $1$ one entry in each of its top row, bottom row, first column, and last column, which we call $t$, $b$, $l$, and $r$ respectively (and these four are distinct).
        \item The $1$ entries $t$ and $b$ are the only $1$ entries in their respective columns, and $l$ and $r$ are the only $1$ entries in their respective rows.
        \item Either $t$ is left of $b$ and $r$ is above $l$, or $t$ is right of $b$ and $r$ is below $l$.
    \end{enumerate}
    Then $P$ has bounded saturation function. (Such patterns are called $Q_1$-like).
\label{thm:q1like}
\end{thm}

This result is enough to completely resolve all patterns with $4$ or fewer ones.

\begin{thm}
    A matrix with four or fewer ones has bounded $\sat(n,P)$ if and only if it is the $1 \times 1$ matrix $\begin{pmatrix}\bullet\end{pmatrix}$ or it is $Q_1$-like.
\end{thm}
\begin{proof}
    
    Assume a matrix $M$ has a bounded saturation function and four or fewer ones. The case where $M$ is $1 \times 1$ is trivial, so we assume $M$ has dimensions greater than $1$. We seek to prove that $M$ cannot contain any $1$ entries in its corners. Without loss of generality, assume that $M$ contains a $1$ in its bottom right corner. By Theorem \ref{thm:decomp}, $M$ is not decomposable. Therefore, $M$ must also contain a second $1$ in its bottom row or rightmost column, because if this were not the case, it would be sum decomposable, which can be seen by letting $B$ equal the $1 \times 1$ matrix $\begin{pmatrix} \bullet \end{pmatrix}$.
    
    Suppose $M$ contains a second $1$ in the bottom row. By Theorem \ref{thm:ssat} there is a $1$ entry in the rightmost column which is the only $1$ entry in its row. Since this cannot be the entry in the corner, there must be a second $1$ entry in the rightmost column.

    Applying Theorem \ref{thm:ssat} to the top row and left column, both must contain a $1$ that is not shared with the bottom row or right column. However, since $M$ has only $4$ ones and we have already used $3$, the only possibility is for there to be a $1$ in the top left corner. However, in this case, the matrix is in fact decomposable, so it has a linear saturation function. 

    Therefore, none of the corners can contain any $1$ entries. This means that all four of the $1$ entries guaranteed by Theorem \ref{thm:ssat} are distinct, so $M$ must be equal to a permutation matrix, possibly with empty rows and columns inserted internally. The only indecomposable permutation matrices (other than $\begin{pmatrix}\bullet\end{pmatrix}$) with four or fewer ones are $Q_1$ and its transpose.
    
    Thus $M$ must be $Q_1$-like; this class of matrices is known to have bounded saturation function by Theorem \ref{thm:q1like}, completing the proof. 
\end{proof}

\subsection{\texorpdfstring{$4 \times 4$}{4 x 4} matrices}

We also complete the classification of the saturation functions of $4 \times 4$ matrices. Note all matrices less than $4 \times 4$ except $\begin{pmatrix} \bullet \end{pmatrix}$ clearly either are decomposable or have linear semisaturation, and hence have linear saturation, so in this sense these are the simplest remaining unresolved cases.

Consider a $4 \times 4$ pattern $P$ with bounded saturation. The $4$ entries on the outside of $P$ guaranteed by Theorem \ref{thm:ssat} are distinct and none lie in the corners; else, we would have a $1$ entry in the corner which is alone in its row and column, and $P$ would be decomposable. In particular, for a $4 \times 4$ pattern $P$, they must lie (up to symmetry) in one of the following two configurations, as observed by Berendsohn in \cite{berendsohn21}.
\begin{align*}
    & \begin{pmatrix}
        \ast&&\bullet&\ast \\
        &&&\bullet \\
        \bullet&&& \\
        \ast&\bullet&&\ast
    \end{pmatrix}, &
    & \begin{pmatrix}
        \ast&&\bullet&\ast \\
        \bullet&&& \\
        &&&\bullet \\
        \ast&\bullet&&\ast
    \end{pmatrix}.
\end{align*}

Since the entry in the first row is the only entry in its column, and similarly the other three, we may conclude that all entries except the corners of $P$ are zeroes. Then only the corners of $P$ are still unknown. At least one of the $1$ entries in $P$ is the only $1$ entry in both its row and column. Clearly this cannot be any of the corners, so it is one of the four which are already marked. The remaining possibilities can each be checked; after eliminating those which are decomposable, we are left with (again up to symmetry) the permutation matrix $Q_1 = \begin{psmallmatrix}
    & & \bullet & \\
    \bullet & & & \\
    & & & \bullet \\
    & \bullet & &
\end{psmallmatrix}$, which is known \cite{berendsohn21} to have bounded saturation, and the following four unresolved cases.
\begin{align*}
    Q_2 &= \begin{pmatrix}
        \bullet & & \bullet & \\
        \bullet & & & \\
        & & & \bullet \\
        & \bullet & &
    \end{pmatrix}, &
    Q_3 &= \begin{pmatrix}
        & \bullet & & \\
        \bullet & & & \\
        & & & \bullet \\
        \bullet & & \bullet & \\
    \end{pmatrix}, \\
    Q_4 &= \begin{pmatrix}
        \bullet & & \bullet & \\
        \bullet & & &  \\
        & & & \bullet \\
        \bullet & \bullet & & \\
    \end{pmatrix}, & Q_5 &= \begin{pmatrix}
         & \bullet & & \\
        \bullet & & & \\
        & & & \bullet \\
        \bullet & & \bullet & \bullet \\
    \end{pmatrix}.
\end{align*}

These are the only four remaining unresolved $4 \times 4$ patterns. We present witnesses for $Q_2$ and $Q_4$, while we show $Q_3$ and $Q_5$ lie in an infinite class of patterns with saturation in $\Theta(n)$. This in particular makes the latter two patterns, along with others in said class, the first known examples of indecomposable matrices with bounded semisaturation function but saturation function in $\Theta(n)$.

\subsection{\texorpdfstring{$Q_2$}{Q2}-like matrices}

Firstly, we resolve the saturation of $Q_2$. In fact, we show further that $Q_2$ lies in an infinite family of matrices with bounded saturation. Say a pattern $P$ is \textit{$Q_2$-like} if
\begin{enumerate}[label=(\roman*)]
    \item $P$ has a $1$ entry in exactly one corner, and
    \item The matrix obtained by deleting this $1$ entry from $P$ is $Q_1$-like.
\end{enumerate}

We present a witness construction. Fix a $Q_2$-like $k \times \ell$ pattern $P$; without loss of generality assume that $P$ has a $1$ in the top left corner. Let the other $1$ entry in the top row of $P$ be $t$, in the first column be $l$, similarly $r$, $b$. Without loss of generality assume $r$ is above $l$, and $t$ is right of $b$. $P$ is of the form
\[P = \begin{pmatrix}
    \bullet & & & & t & & & \\
    & A & & B & & C & \\
    l & & & & & & \\
    & D & & E & & F & \\
    & & & & & & r \\
    & G & & H & & I & \\
    & & b & & & &
\end{pmatrix},\]
for any matrices $A,\dots,I$ of suitable dimensions (i.e.\ $A,D,G$ have the same number of columns).

Let $h_l$ be the row containing the entry $l$, $c_t$ the column containing $t$, similary $h_r$, $c_b$. Denote by $P_t$, $P_l$, $P_r$, and $P_b$ the patterns obtained by deleting the top, left, right, and bottom rows of $P$ respectively.

Recall from \cite{berendsohn21} the construction $W_1(P)$, which consists of a copy of $P_l$ and $P_r$ offset from each other. Define $W_V(P)$ as follows: we begin with $W_1(P)$, and then, at the row of $W_1(P)$ corresponding with the top of $P_l$ (row $h_l - h_r$), for every entry of this row that lies in the copy of $P_r$ (i.e. in the first $\ell - 1$ columns of this row), if there is no $1$ entry, then we add a $1$. To illustrate,

\[W_V(P) =
\begin{pmatrix}
    \bullet & & & & t & & \\
    & A & & B & & C \\
    \bullet & \cdots & \bullet & \cdots & \bullet & \cdots & & & & t & & \\
    & D & & E & & F & A & & B & & C \\
    \cdot & \cdot & \cdot & \cdot & \cdot & \cdot & \cdot & \cdot & \cdot & \cdot & \cdot & \cdot\\
    & G & & H & & I & D & & E & & F\\
    & & b & & & & & & & & & r\\
    & & & & & & G & & H & & I \\
    & & & & & & & b & & &\\
\end{pmatrix},
\]
noting that we slightly abuse notation here, as $A,B,C$ may have different numbers of rows than $D,E,F$ and $G,H,I$, so, for example, while $b$ and $r$ are shown aligned in the above illustration, they may not be in general.

Now define $W_H(P)$ as follows. We begin by placing a copy of $P_b$ above a copy of $P_t$, such that the columns $c_b$ and $c_t$ of each respectively are aligned, similarly to the construction in \cite{berendsohn21} of $W_1$. We remove the first column of this pattern (which is empty except for the entry $l$ in the copy of $P_t$); this is $W_{H0}(P)$. 
%
Note that $W_{H0}(P)$ has $2k - 2$ rows, the top $k - 1$ of which correspond to $P_b$ and the bottom $k - 1$ of which correspond to $P_t$. Now, using block matrix notation, define the $(2k - 2) \times (k - 1)$ pattern $X(P)$ as follows, where, as above, $\mathbf{1}_{m\times n}$ is the $m \times n$ all-ones matrix.
\[X(P) =
\begin{pmatrix}
    I_{k - 1} \\ \mathbf{0}_{(h_l - 2) \times (k - 1)} \\ \mathbf{1}_{1\times(k-1)} \\ \mathbf{0}_{(k - h_l) \times (k - 1)}\\ 
\end{pmatrix}.\]
Then, we let \begin{align*}
    W_H(P) &=
    \begin{pmatrix}
        X(P) & W_{H0}(P)
    \end{pmatrix} \\
    &= \begin{pmatrix}
        \bullet & & & & & & & \bullet & & \cdot & & t & & & \\
        & \ddots & & & & & & & A & \cdot & B & & C & \\
        & & \bullet & & & & & l & & \cdot & & & & \\
        & & & \ddots & & & & & D & \cdot & E & & F & \\
        & & & & \bullet & & & & & \cdot & & & & r \\
        & & & & & \ddots & & & G & \cdot  & H & & I & \\
        & & & & & & A & & B & \cdot  & C & & &\\
        \bullet & \cdots & \bullet & \cdots & \bullet & \cdots & & & & \cdot & & & & \\
        & & & & & & D & & E & \cdot & F & & & \\
        & & & & & & & & & \cdot & & r & & \\
        & & & & & & G & & H & \cdot & I & & & \\
        & & & & & & & b & & \cdot & & & &
    \end{pmatrix}.
\end{align*}

(As above, we abuse notation; $t$ and $r$ need not be aligned as they are in the diagram, nor must $l$ and $b$). Note that the row of all ones in the lower half of $X(P)$ is aligned with the row in $W_{H0}$ of the copy of $P_t$ which would contain the entry $l$.
\begin{example}
    For $Q_2$, the simplest $Q_2$-like pattern, we have the witnesses
    \[W_V(Q_2) = \begin{pmatrix}
        \bullet & & \bullet & & \\
        \bullet & \bullet & \bullet & & \bullet & \\
        \cdot & \cdot & \cdot & \cdot & \cdot & \cdot \\
        & \bullet & & & & \bullet \\
        & & & \bullet & &
    \end{pmatrix}, W_H(Q_2) = \begin{pmatrix}
        \bullet & & & \bullet & \cdot & \bullet & \\
        & \bullet & & \bullet & \cdot & & \\
        & & \bullet & & \cdot & & \bullet \\
        \bullet & \bullet & \bullet & & \cdot & \\
        & & & & \cdot & \bullet \\
        & & & \bullet & \cdot \\
    \end{pmatrix}.\]
\end{example}

It may be verified manually (or using the code in \cite{code}) that these particular examples are witnesses. Indeed the construction gives witnesses in general, as we show.

\begin{theorem} \label{thm:q2like}
    If $P$ is $Q_2$-like, then $\sat(n,P) = O(1)$.
\end{theorem}

\begin{proof}
    By Theorem \ref{thm:vhwitness}, it suffices to present vertical and horizontal witnesses. We claim that $W_V(P)$ is a vertical witness for $P$, and $W_H(P)$ is a horizontal witness.

    Firstly, we consider $W_V(P)$. We will see that the expandable row is row $h_r$, i.e. the empty row which was aligned in the copies of $P_r$ and $P_l$ in the construction. Fix any column $1 \leq c \leq 2\ell - 2$ of $W_V(P)$, and change the entry in row $h_r$ of $c$ to a $1$. If $c \geq \ell$, i.e. $c$ is in the right half of $W_V(P)$, then $c$ together with $P_r$ forms a copy of $P$. If $c$ is in the left half, $c \leq \ell - 1$, then note there is a $1$ entry in column $c$ in row $h_r - h_l$, the top row of the copy of $P_l$. Then column $c$, together with $P_l$, must form a copy of $P$ (noting that the $1$ entries in $c$ correspond to the corner of $P$ and entry $l$ respectively). Thus row $h_r$ is expandable.

    It remains to see that the pattern $W_V$ does not contain $P$. Suppose otherwise, that $W_V$ contains some copy $P_0$ of $P$, with inclusion map $\phi$. Then $\phi(t)$ must have at least $k - 1$ nonempty rows below it. Note that the copy of $P_l$ has only $k - 1$ nonempty rows total, so $\phi(t)$ cannot be in the copy of $P_l$. Then $\phi(t)$ must be in the left half of $W_V$, the copy of $P_r$. Then $\phi(b)$ is left of $\phi(t)$, so $\phi(b)$ is also in the left half of $W_V$. But also $\phi(b)$ is at or below row $k$, but there are no nonempty entries of $W_V$ in the left half below row $k$, a contradiction; thus $W_V$ does not contain $P$, so $W_V$ is indeed a vertical witness for $P$.

    Now, we consider $W_H(P)$. We first show that column $\ell + c_t - 1$ (the column marked with dots in the diagram) is expandable. It follows from the construction of $W_{H0}$ that this column is empty. Suppose we change the entry in this column in row $r \in [1,2k - 2]$ to a $1$. If $k \leq r$, then we can observe that $r$, together with the copy of $P_b$, form a copy of $P$. This leaves the case $r \leq k - 1$. In this case, we consider the $k$ rows $r$ and $k,\dots,2k - 2$, and the $\ell$ columns $r,k,\dots,k+\ell-1$. Note that these are the rows and columns corresponding to the copy of $P_t$, together with row $r$ and column $r$. Note that these additional rows and columns have entries corresponding exactly to $l$, $t$, and the entry in the corner, so we see that this is indeed a copy of $P$.

    Lastly, we must show that $W_H(P)$ does not contain $P$. Suppose that $W_H(P)$ contains a copy $P_0$ of $P$ with inclusion map $\psi$. Then $\psi(l)$ must have at least $\ell - 1$ columns to its left, and the column $\psi(c_b)$ cannot be empty, together meaning $\psi(l)$ cannot be in the copy of $P_b$. We also know $\psi(l)$ has a $1$ entry above it, so it cannot be in the copy of $I_{k-1}$ in $X(P)$. Together, this means $\psi(l)$ cannot be in the top $k - 1$ rows. Now, $\psi(r)$ is below $\psi(l)$, so it must also be in the bottom $k - 1$ rows. Then $\psi(r)$ cannot be right of column $k + \ell - 1$ (the rightmost column of the copy of $P_t$). Also, $\psi(r)$ must be at least $\ell - 1$ rows to the right of $\psi(l)$, meaning $\psi(l)$ must be left of column $k$. Thus, $\psi(l)$ is in $X(P)$, and in particular in the row of all ones in $X(P)$ (since we previously eliminated $I_{k-1}$). Then $\psi(b)$ is at least $k - h_l$ rows below this, meaning $\psi(b)$ must be the bottom entry of the pattern. Then $\psi(r)$ is at least $\ell - c_b$ columns to the right of this; then $\psi(r)$ can only be in column $\ell + k - 1$, and every column between $\psi(r)$ and $\psi(b)$ must be in $P_0$. But then $\psi(t)$ must be in column $\ell + c_t - 1$, which is empty - a contradiction. Therefore, $W_H(P)$ does not contain $P$, so $W_H(P)$ is a horizontal witness for $P$.

    Thus, $\sat(n,P) = O(1)$.
\end{proof}

\begin{figure}
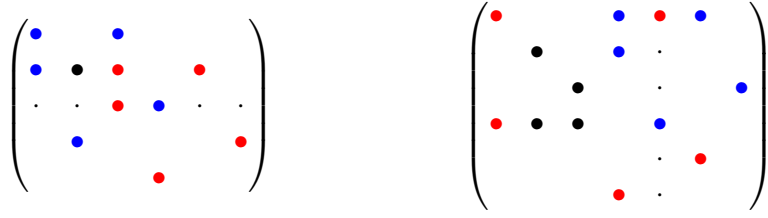

    \begin{align*}
        &\begin{pmatrix}
            \textcolor{blue}{\bullet} & & \textcolor{blue}{\bullet} & & \\
            \textcolor{blue}{\bullet} & \bullet & \textcolor{red}{\bullet} & & \textcolor{red}{\bullet} & \\
            \cdot & \cdot & \textcolor{red}{\bullet} & \textcolor{blue}{\bullet} & \cdot & \cdot \\
            & \textcolor{blue}{\bullet} & & & & \textcolor{red}{\bullet} \\
            & & & \textcolor{red}{\bullet} & &
        \end{pmatrix} & 
        &\begin{pmatrix}
            \textcolor{red}{\bullet} & & & \textcolor{blue}{\bullet} & \textcolor{red}{\bullet} & \textcolor{blue}{\bullet} & \\
            & \bullet & &\textcolor{blue}{\bullet} & \cdot & & \\
            & & \bullet & & \cdot & & \textcolor{blue}{\bullet} \\
            \textcolor{red}{\bullet} & \bullet & \bullet & & \textcolor{blue}{\bullet} & & \\
            & & & & \cdot & \textcolor{red}{\bullet} & \\
            & & & \textcolor{red}{\bullet} & \cdot & &
        \end{pmatrix}
    \end{align*}
    \caption{Examples of $Q_2$ in $W_V(Q_2)$ and $W_H(Q_2)$}
    \label{fig:witnesses}
\end{figure}

\begin{prop}
    The saturation function of $Q_4$ is in $O(1)$.
\end{prop}

\begin{proof}
\label{prop:q4}
    We first show a vertical and horizontal witness for $Q_4$. Let
    \[W_V = .\begin{pmatrix}
         \bullet&&\bullet&&&\\
        \bullet&\bullet&\bullet&&\bullet&\\
        \cdot&\cdot&\cdot&\cdot&\cdot&\cdot\\
        \bullet&\bullet&&&&\bullet\\
        &\bullet&\bullet&\bullet&&\\
        \bullet&&\bullet&&&\\
    \end{pmatrix}
    W_H = \begin{pmatrix}
        \bullet&&&\bullet&\cdot&\bullet&\\
        &\bullet&&\bullet&\cdot&&\\
        \bullet&\bullet&&&\cdot&&\bullet\\
        &&&&\cdot&&\\
        &\bullet&\bullet&\bullet&\cdot&&\\
        &\bullet&&\bullet&\cdot&\bullet&\\
        &\bullet&&\bullet&\cdot&&\\
        \bullet&&&\bullet&\cdot&&
    \end{pmatrix}
    \]
    We then glue these together using the construction described in Theorem \ref{thm:str_indecomp}; though $Q_4$ is not strongly indecomposable, we will observe that the construction still gives a valid witness. In particular, we verify using our code \cite{code} that the matrix
    \[W = \begin{pmatrix}
        &&&&&&\bullet&&&\bullet&\cdot&\bullet&\\
        &&&&&&&\bullet&&\bullet&\cdot&&\\
        &&&&&&\bullet&\bullet&&&\cdot&&\bullet\\
        \bullet&&\bullet&&&&&&&&\cdot&&\\
        \bullet&\bullet&\bullet&&\bullet&&&&&&\cdot&&\\
        \cdot&\cdot&\cdot&\cdot&\cdot&\cdot&\cdot&\cdot&\cdot&\cdot&\cdot&\cdot&\cdot\\
        \bullet&\bullet&&&&\bullet&&&&&\cdot&&\\
        &\bullet&\bullet&\bullet&&&&&&&\cdot&&\\
        \bullet&&\bullet&&&&&&&&\cdot&&\\
        &&&&&&&\bullet&\bullet&\bullet&\cdot&&\\
        &&&&&&&\bullet&&\bullet&\cdot&\bullet&\\
        &&&&&&&\bullet&&\bullet&\cdot&&\\
        &&&&&&\bullet&&&\bullet&\cdot&&
    \end{pmatrix}\] is a witness for $Q_4$.
\end{proof}

\subsection{\texorpdfstring{$Q_3$}{Q3}-like matrices}

We now present an infinite class of patterns with linear saturation function - in particular, a class including $Q_3$ and $Q_5$, among other patterns with bounded semisaturation function and which are indecomposable. In particular, call a pattern $P$ \textit{$Q_3$-like} if it is of the form
\[P = \begin{pmatrix}
    0 & 1 & \mathbf{0} & 0 \\
    1 & 0 & \mathbf{0} & 0 \\
    0 & 0 & \mathbf{0} & 1 \\
    A & \mathbf{0} & B & C
\end{pmatrix}.\] for some submatrices $A$, $B$, and $C$ of appropriate dimensions, where $\mathbf{0}$ represents a $0$ matrix of any size, where we also require require at least two $1$ entries in all rows but the first three. In other words, a $k \times \ell$ pattern $P$ is $Q_3$-like if
    \begin{enumerate}[label=(\roman*)]
        \item The entries $p_{12} = p_{21} = p_{3\ell} = 1$ are all $1$ entries and are the only $1$ entries in their rows,
        \item $p_{12}$ is the only $1$ entry in the second column, and
        \item all rows of $P$ other than the first $3$ have at least two $1$ entries.
    \end{enumerate}
We have the following result for $Q_3$-like patterns.

\begin{prop}
\label{prop:q3_like}
    Let $P$ be any pattern. If $P$ is $Q_3$-like, then $\sat(m,n_0,P) = \Theta(n)$.
\end{prop}

\begin{proof}
    Let the $k \times \ell$ pattern $P$ be $Q_3$-like. Suppose that $\sat(m,n_0,P)$ is bounded; then in particular we must have a vertical witness $W = \begin{pmatrix} w_{i,j} \end{pmatrix}$ for $P$. Since $P$ has no empty rows, we may assume without loss of generality that the only empty row of $W$ is a single expandable row $i$ (else, the matrix formed by deleting extraneous empty rows of $W$ clearly remains a witness for $P$).
    
    \begin{figure}
        \centering
        $\begin{pmatrix}
            & & & & \textcolor{red}{\circ} & & & & \\
            & \textcolor{red}{\bullet} & & \circ & \textcolor{blue}{\bullet} & & \bullet & & \\
            \cdot & w_{i,j_3} & \cdot & w_{i,a_1} & \cdot & \cdot & w_{i,j_2} & \cdot \\
            & \textcolor{red}{\bullet} & & & & \textcolor{red}{\bullet} & \textcolor{blue}{\circ} & & \\
            & & & \textcolor{blue}{\circ} & & \textcolor{blue}{\circ} & & & \\
        \end{pmatrix}$
        \caption{Finding $P_0$ in an example for $Q_3$ as described in Proposition \ref{prop:q3_like}. The one entries in red are those in $P_{j_2}$, while those in blue are those in $P_{a_1}$. The entries marked by $\circ$ are those in $P_0$.}
        \label{fig:q3_construction}
    \end{figure}
    
    We denote $t, l, r$ for the entries $p_{12}$, $p_{21}$, and $p_{3\ell}$ of $P$ respectively. For any $j$, as $w_{i,j}$ is the only $1$ entry in its row of $P_j$, it must be in one of the top $3$ rows of $P_j$; thus $w_{i,j}$ is $\phi_j(t)$, $\phi_j(r)$, or $\phi_j(l)$.

    If $i = 1$, then $P_1$ would have a $1$ entry in its top left corner, which is a contradiction. Thus we have a row $i - 1$. Furthermore, row $i - 1$ is nonempty by our construction of $W$, so we define the set of columns $J = \{j : w_{i + 1, j} = 1\} \neq \emptyset$. For any $j \in J$, if $w_{i,j} = \phi_j(t)$, then we could replace the one entry $w_{i,j}$ by $w_{i-1,j}$ in $P_j$ (i.e., replace the row $i$ by $i - 1$) and obtain a copy of $P$ in $W$ itself, which is a contradiction. Thus we must have either $w_{i,j} = \phi_j(r)$ or $\phi_j(l)$. We say the column $j$ \textit{corresponds} to either $r$ or $l$.
    
    Since $J$ is nonempty, let $j_0$ be the rightmost column in $J$. Suppose $j_0$ corresponds to $l$. Let $w_{i_1, j_1} = \phi_{j_0}(t)$; clearly $i_1 < i$, $j_1 > j_0$. If $i_1 < i - 1$, then replacing row $i$ by row $i - 1$ in $P_{j_0}$ gives a copy of $P$ in $W$, which is a contradiction. Then we must have $i_1 = i - 1$; that is, we have a $1$ entry $w_{i - 1, j_1} = 1$ in $W$, so $j_1 \in J$. But $j_1 > j_0$, which contradicts the fact the $j_0$ was the rightmost column in $J$. Thus in fact $j_0$ corresponds to $r$.

    Thus we have at least one column in $J$ corresponding to $r$. Let $j_2$ be the leftmost such column; then $w_{i,j_2}$ is contained in a copy $P_{j_2}$ of $P$ in $W_{j_2}$, and in particular $w_{i,j_2} = \phi_{j_2}(r)$. Let $w_{i_3,j_3} = \phi_{j_2}(l)$ and $w_{i_4,j_4} = \phi_{j_2}(t)$. We have $i_4 < i_3 < i$ and $j_2 < j_4 < j_3$. Suppose that $i_3 < i - 1$; then, because $w_{i - 1,j_2} = 1$, replacing the row $i$ by $i - 1$ in $P_{j_2}$ gives a copy of $P$ in $W$, which is a contradiction, so we must have in fact $i_3 = i - 1$. 
    
    Then $j_3 \in J$, so we have at least one column in $J$ left of $j_4$. Let $a_1$ be the rightmost column in $J$ left of $j_4$. As $j_2$ was the leftmost column corresponding to $r$ and $a_1 < j_2$, the column $a_1$ must correspond to $l$. Let the $\ell$ columns of $P_{a_0}$ be $a_1 < a_2 < \cdots < a_\ell$, and let the $k$ rows be $b_1 < b_2 = i < \cdots < b_k$. Suppose $b_1 < i - 1$; then, similarly to above, because $w_{i-1,a_1} = 1$, replacing row $i$ of $P_{a_1}$ by $i - 1$ gives a copy of $P$ in $W$, which is a contradiction. Thus $b_1 = i - 1$, so in particular since $w_{b_1,a_2} = 1$, we have $a_2 \in J$. Then since $a_1$ was the rightmost column in $J$ left of $j_4$, we must have $a_\ell > \cdots > a_2 \geq j_4$.

    Consider the submatrix $P_0$ (and associated embedding $\phi_0$) of $W$ formed by the columns $a_1 < j_4 < a_3 < \cdots < a_\ell$ and the rows $i_4 < i - 1 < b_3 < \cdots < b_k$. We claim that $P_0$ is a copy of $P$ in $W$. Since $P_{a_1}$ is a copy of $P$, for any $\iota \notin \{1,2\}$ and $\jmath \neq 2$, we know $\phi_0(p_{\iota,\jmath}) = \phi_{a_1}(p_{\iota,\jmath})$ is a $1$ entry. It remains to check the second column and the first two rows. By our assumption on $P$, these are empty apart from the two one entries $t$ and $l$. But we have seen $\phi_0(t) = w_{i_4,j_4} = \phi_{j_2}(t) = 1$ and $\phi_0(l) = w_{i - 1, a_1}$, which, because $a_1 \in J$ by assumption, is also a $1$ entry.
    
    Therefore, $P_0$ is a copy of $P$ in $W$, which is a contradiction. Thus, $\sat(n,P) = \Theta(n)$.
\end{proof}

Since $Q_3$ and $Q_5$ are both $Q_3$-like, this completes the classification of $4 \times 4$ patterns by boundedness of $\sat(n,\cdot)$. However, many patterns still have unresolved $\sat(m_0,n,\cdot)$ or $\sat(m,n_0,\cdot)$. In particular, some patterns have both $\sat(m_0,n,\cdot)$ and $\sat(m,n_0,\cdot)$ bounded, but $\sat(n,\cdot)$ linear. A wrinkle in this classification is illustrated by the following example.

\begin{example}
\label{ex:q6}
    The patterns \begin{align*}
    \begin{pmatrix}
        \bullet & & \bullet & \\
        \bullet & & & \\
        & & & \bullet \\
        \bullet & \bullet & & \bullet
    \end{pmatrix} & &\begin{pmatrix}
        \bullet & & \bullet & \\
        \bullet & & & \\
        & & & \bullet \\
        & \bullet & & \bullet
    \end{pmatrix} & &\begin{pmatrix}
        \bullet & & \bullet & \\
        & & & \bullet \\
        \bullet & & & \\
        & \bullet & & \bullet
    \end{pmatrix}
\end{align*} have linear $\sat(n, \cdot)$ by Theorem~\ref{thm:ssat}, despite also having both vertical and horizontal witnesses (see Table~\ref{tab:witnesses}).
\end{example}

We also present witnesses for a new family of matrices.

\begin{obs}\label{obs:intermediaries}
    For any patterns $P$, $Q$, where $P$ is contained in $Q$, suppose that $W$ is a matrix avoiding $P$ with a $Q$-expandable row (column). Then for every intermediary pattern $R$ contained in $Q$ and containing $P$, $W$ is a vertical (horizontal) witness for $R$.
\end{obs}

This gives a new method of constructing witnesses for entire families of matrices. For example, we can use this to give witnesses for $3$ new $5\times5$ matrices (up to symmetry).

\begin{example}
    Let
    \[Q_9 = \begin{pmatrix}
        \bullet&\bullet& & &\bullet\\
         & & & &\bullet\\
         & &\bullet& & \\
        \bullet& & & & \\
        \bullet& & &\bullet&\bullet
    \end{pmatrix},\hspace{1cm}
    P = \begin{pmatrix}
        &\bullet& & &\bullet\\
         & & & &\bullet\\
         & &\bullet& & \\
        \bullet& & & & \\
        \bullet& & &\bullet&
    \end{pmatrix}.\]
    Since
    \[W_V=\begin{pmatrix}
        \bullet&\bullet&\bullet&&&\bullet&&&&&&\\
        &&&&&\bullet&&&&&&\\
        \bullet&\bullet&&&&&&&&&&\bullet\\
        &&&\bullet&&&\bullet&\bullet&&&\bullet&\bullet\\
        \cdot&\cdot&\cdot&\cdot&\cdot&\cdot&\cdot&\cdot&\cdot&\cdot&\cdot&\cdot\\
        \bullet&\bullet&&&\bullet&\bullet&&&\bullet&&&\\
        \bullet&&&&&&&&&&\bullet&\bullet\\
        &&&&&&\bullet&&&&&\\
        &&&&&&\bullet&&&\bullet&\bullet&\bullet
    \end{pmatrix}\]
    has a $Q_9$-expandable row and avoids $P$ \cite{code}, we know that any matrix of the form
    \[R = \begin{pmatrix}
        \ast&\bullet& & &\bullet\\
         & & & &\bullet\\
         & &\bullet& & \\
        \bullet& & & & \\
        \bullet& & &\bullet&\ast
    \end{pmatrix}\]
    has $W_V$ as a vertical witness. Therefore, $\sat(n,m_0,R)=O(1)$. 
\end{example}

\begin{remark}
    If we rotate $W_V$ by $90^\circ$, then we also obtain a horizontal witness for $Q_9$ (but it contains $P$ unfortunately). Since the entries in the expandable row of $W_V$ that are at least two entries from the edge of the matrix all correspond to the central $1$ entry in $Q_9$, we can use Theorem \ref{thm:str_indecomp} to deduce that $\sat(n,Q_9)=O(1)$.
\end{remark}

As another application of Observation \ref{obs:intermediaries}, we show that $k \times 4$ matrices of the form

\[R_k = \begin{pmatrix}
    \bullet & & \bullet & \bullet \\
    \ast & & & \ast \\
    \vdots & & & \vdots \\
    \ast & & & \ast \\
    \bullet & \bullet & & \bullet
\end{pmatrix}\] have bounded saturation. These are intermediaries between\[P_k = \begin{pmatrix}
    \bullet & & \bullet & \bullet \\
     & & & \\
    \vdots & & & \vdots \\
     & & & \\
    \bullet & \bullet & & \bullet
\end{pmatrix} \text{ and } Q_k = \begin{pmatrix}
    \bullet & & \bullet & \bullet \\
    \bullet & & & \bullet \\
    \vdots & & & \vdots \\
    \bullet & & & \bullet \\
    \bullet & \bullet & & \bullet
\end{pmatrix}.\] Define \[W_2 = \begin{pmatrix}
    & \bullet & \bullet & \cdot & \bullet & \bullet \\
    \bullet & \bullet & \bullet & \cdot & \bullet & \\
    & & \bullet & \cdot & \bullet & \bullet \\
    \bullet & & \bullet & \cdot & \bullet & \bullet
\end{pmatrix}.\]
Then, for any $k > 2$, define $W_k = W_2 \otimes \mathbf{1}_{k-1}$, where $\mathbf{1}_{k-1}$ is the all-ones column vector with $k - 1$ entries. We will see that $W_k$ satisfies the conditions of Observation \ref{obs:intermediaries}.

\begin{prop}\label{prop:intermediaries_1}
    Any $k \times 4$ pattern $R_k$ as defined above has $\sat(m_0,n,R_k) = O(1)$.
\end{prop}


\begin{proof}
    The case $k = 2$ is fairly simple, as only one possible pattern $R_2$ exists to check: $R_2 = P_2 = Q_2$. It is clear that $W_2$ does not contain $R_2$.
    
    Now, we check in detail that the fourth column of $W_2$ is expandable for $R_2$. In particular, if we change the entry $o_{1,4}$ in row $1$ and column $4$ to a 1, then there is clearly a copy of $Q_2$ in rows $1$ and $2$, and columns $2$, $3$, $4$, $5$. If we change the entry $o_{2,4}$ to a $1$, then we introduce a copy of $Q_2$ in rows $2$ and $4$ and columns $1$, $3$, $4$, $5$; changing $o_{3,4}$ gives a copy of $Q_2$ in rows $1$ and $3$ and columns $3$, $4$, $5$, $6$; and changing the fourth entry in row $4$ to a $1$ gives a copy of $Q_2$ in rows $3$ and $4$, and columns $3$, $4$, $5$, $6$.
    
    Let $C_i$ be the set of columns in the copy of $Q_2$ introduced by changing the entry $o_{i,4}$ to a $1$, so $C_1 = \{2,3,4,5\}$, and so on; we will use these in the general case.

    For the case $k > 2$, we must show that $W_k$ avoids $P_k$ and has a $Q_k$-expandable column.
    
    We first show the latter. Let $i$ be any row of $W_k$, and create the matrix $(W_k)_i$ by changing the entry of column $4$ in row $i$ to a $1$. By construction of $W_k$, row $i$ of $W_k$ is equal to one of the four rows of $W_2$; in particular say it is equal to row $i'$ of $W_2$. Per our above observations in the case of $W_2$, there is a $2 \times 4$ submatrix $(Q_2)_i$ of $(W_2)_i$ corresponding to a copy of $Q_2$.
    
    The submatrix $(Q_2)_i$ has two rows: row $i'$ and one other row, which we call row $j' \in \{1,2,3,4\}$. There are $k - 1$ consecutive rows of $W_k$ which are equal to row $j'$ of $W_2$. Let the submatrix $(Q_k)_i$ of $(W_k)_i$ contain these $k - 1$ rows as well as row $i$, and the columns given by $C_i$ (that is, the same columns as $(Q_2)_i$). In particular we will label the first column $a = \min C_i$ and the last column $d = \max C_i$.
    
    We must check that $(Q_k)_i$ is a copy of $Q_k$ in $(W_k)_i$. The top and bottom rows of $(Q_k)_i$ are exactly a copy of $Q_2$, as expected; one is row $i$ and the other is a copy of row $j'$. Consider any row $j$ of $(Q_k)_i$ between the top and bottom; $j$ is equal to the row $j'$ of $(Q_2)_{i'}$. This corresponds to either the first or second row of $Q_2$; in either case, both the first and last entries of row $j'$ in $(Q_2)_{i'}$ are one entries. That is, the entries of row $j$ of $(W_k)_i$ in columns $a$ and $d$ are both one entries, so $(Q_k)_i$ is indeed a copy of $Q_k$.
    
    Lastly, suppose $W_k$ contains a copy $P_k^0$ of $P_k$. Each row of $W_k$ is equal to one of the four rows of $W_2$; in particular, as there are at least $k - 1$ rows between the first and last row of $P_k^0$, they must correspond to different rows of $W_2$. Then these two rows in $W_2$ would in fact contain a copy of $P_2$, but it is easy to verify that $W_2$ avoids $P_2$ \cite{code}, which completes the proof.
\end{proof}

To finish, we give a table of various other miscellaneous patterns with witnesses for them. In each row, $W_V$ is a vertical witness for $P$ and $W_H$ is a horizontal witness for $P$. When $W_H$ is marked $W_V^T$, the pattern is symmetric and so we can obtain a horizontal witness simply by transposing the vertical witness. Each of these witnesses has been verified with the code in \cite{code}.

\newpage

\begin{table}
\centering
\begin{tabular}{c|c|c}
    $P$ & $W_V$ & $W_H$ \\ \midrule
    $\begin{psmallmatrix}
        \bullet & & \bullet & \\
        \bullet & & & \\
        & & & \bullet \\
        \bullet & \bullet & & \bullet
    \end{psmallmatrix}$ & $\begin{psmallmatrix}
        \bullet&&\bullet&&&&\\
        \bullet&\bullet&\bullet&&\bullet&&\\
        \cdot&\cdot&\cdot&\cdot&\cdot&\cdot&\cdot\\
        \bullet&\bullet&&\bullet&\bullet&\bullet&\bullet\\
        &\bullet&\bullet&\bullet&&\bullet&\\
        \bullet&&\bullet&&&&\bullet\\
    \end{psmallmatrix}$ & $\begin{psmallmatrix}
        \bullet & & \bullet & \cdot & \bullet & \\
        & \bullet & \bullet & \cdot & & & \\
        & \bullet & & \cdot & & \bullet & \\
        & \bullet & \bullet & \cdot & & \bullet \\
        \bullet & & \bullet & \cdot & \bullet & \bullet & \bullet \\
        & \bullet & \bullet & \cdot & \bullet & \bullet & \\
        \bullet & & \bullet & \cdot & \bullet & & \bullet 
    \end{psmallmatrix}$ \\ \midrule
    
    $\begin{psmallmatrix}
        \bullet & & \bullet & \\
        \bullet & & & \\
        & & & \bullet \\
        & \bullet & & \bullet
    \end{psmallmatrix}$ & $\begin{psmallmatrix}
        \bullet & & \bullet & & & \\
        \bullet & \bullet & \bullet & & \bullet & \\
        \cdot & \cdot & \cdot & \cdot & \cdot & \cdot \\
        & \bullet & & \bullet & \bullet & \bullet \\
        & & & \bullet & & \bullet
    \end{psmallmatrix}$ & $\begin{psmallmatrix}
        \bullet & & \bullet & \cdot & \bullet & & \\
        & \bullet & \bullet & \cdot & & & \\
        & \bullet & & \cdot & & \bullet & \bullet \\
        \bullet & \bullet & & \cdot & & \bullet & \\
        & & & \cdot & \bullet & \bullet & \\
        & & \bullet & \cdot & \bullet & & \bullet
    \end{psmallmatrix}$ \\ \midrule
    
    $\begin{psmallmatrix}
        \bullet & & \bullet & \\
        & & & \bullet \\
        \bullet & & & \\
        & \bullet & & \bullet
    \end{psmallmatrix}$ & $\begin{psmallmatrix}
        & & & & & \bullet & \bullet & & \bullet & \\
        & & \bullet & & \bullet & & & & & \\
        & \bullet & & \bullet & & & & & & \\
        \bullet & & \bullet & & & & & & & \\
        \bullet & & \bullet & & & & \bullet & \bullet & \bullet & \bullet \\
        \cdot & \cdot & \cdot & \cdot & \cdot & \cdot & \cdot & \cdot & \cdot & \cdot \\
        \bullet & \bullet & \bullet & \bullet & & & & \bullet & & \bullet \\
        & & & & & & & \bullet & & \bullet \\
        & & & & & & \bullet & & \bullet & \\
        & & & & & \bullet & & \bullet & & \\
        & \bullet & & \bullet & \bullet & & & & &
    \end{psmallmatrix}$ & $W_V^T$ \\  \midrule
    
    $\begin{psmallmatrix}
        \bullet & & \bullet & \bullet \\
        & & & \bullet \\
        \bullet & & & \\
        \bullet & \bullet & & 
    \end{psmallmatrix}$ & $\begin{psmallmatrix}
        & \bullet & \bullet & \bullet & \bullet & \bullet & \bullet \\
        \bullet & & \bullet & \bullet & & & \\
        \bullet & & & \bullet & \bullet & \bullet & \bullet \\
        \cdot & \cdot & \cdot & \cdot & \cdot & \cdot & \cdot \\
        \bullet & \bullet & \bullet & \bullet & \bullet & & \\
        \bullet & & \bullet & & & &
    \end{psmallmatrix}$ & $W_V^T$ \\ \midrule
    $\begin{psmallmatrix}
        \bullet & \bullet & & \bullet \\
        \bullet & & & \\
        & & & \bullet \\
        \bullet & & \bullet &
    \end{psmallmatrix}$ & $\begin{psmallmatrix}
        & \bullet & \bullet & & \bullet & & \\
        & & \bullet & \bullet & & \bullet & \\
        \bullet & & & & \bullet & & \bullet \\
        \bullet & \bullet & \bullet & \bullet & \bullet & & \bullet \\
        \cdot & \cdot & \cdot & \cdot & \cdot & \cdot & \cdot \\
        \bullet & \bullet & \bullet & \bullet & \bullet & \bullet & \bullet \\
        & & & \bullet & & \bullet & \\
        & & \bullet & & & \bullet & \\
        & \bullet & & & & \bullet & \\
        \bullet & & & \bullet & & &
    \end{psmallmatrix}$ & $W_V^T$ \\ \midrule
    $\begin{psmallmatrix}
        \bullet & & \bullet & \\
        \bullet & & & \bullet \\
        & \bullet & & \bullet
    \end{psmallmatrix}$ & - & $\begin{psmallmatrix}
        & & \bullet & \cdot & \bullet & \bullet & \\
        \bullet & \bullet & & \cdot & \bullet & & \bullet \\
        \bullet & & \bullet & \cdot & & \bullet & \bullet \\
        & \bullet & \bullet & \cdot & \bullet & &
    \end{psmallmatrix}$
\end{tabular}
\caption{Various witnesses of matrices that have a linear saturation in two dimensions but not in one dimension.}
\label{tab:witnesses}
\end{table}

It seems likely that many of these may be generalizable to infinite families in the same vein as Theorem~\ref{thm:q2like} and the like.

\section{Discussion}\label{sec:conc}

In this paper, we proved many results about saturation and semisaturation for 0-1 matrices. In particular, we determined the effects of a number of basic operations on the saturation and semisaturation functions, including adding empty columns and rows, and performing Kronecker products. We showed a number of results about saturation functions and semisaturation functions with one fixed dimension, extending several results from the case where both dimensions have the same length. We introduced witness graphs, which we used to determine the saturation functions of some previously unresolved patterns. We also found a linear program for computing the saturation function, and we used it to compute the saturation function for small values of $n$ and make some more general conjectures. Finally, we determined the saturation functions of many small, previously unresolved patterns up to a constant factor. In particular, we completed the classification of the saturation functions of $4 \times 4$ matrices, along with many others.

In Theorem~\ref{thm:str_indecomp}, we characterized the strongly indecomposable matrices with bounded saturation functions. However, there still exist matrices that are not strongly indecomposable for which it is possible to construct a witness using the horizontal and vertical witnesses in a similar fashion to the proof of Theorem~\ref{thm:str_indecomp}, such as for $Q_4$ in Proposition~\ref{prop:q4}. The problem of categorizing such non-strongly indecomposable matrices remains open. 

Additionally, while we have proven by construction in Example~\ref{ex:q6} that there exist matrices with linear saturation function despite also having vertical and horizontal witnesses, the manner by which it was determined that these matrices do indeed have linear saturation functions was by leveraging the fact that they have linear semisaturation functions. This prompts the natural question of whether or not there exists a matrix that has a constant semisaturation function and vertical and horizontal witnesses, but linear saturation function. In other words, can the condition of strong indecomposability and the conditions on the witnesses in Theorem \ref{thm:str_indecomp} be weakened to imply that the semisaturation function is bounded? This question also remains open.

While Theorem~\ref{thm:3d_bounded_sat} gives a $3$-dimensional pattern with $\sat(m_1,m_2,n; P) = O(1)$, there are still many questions for multidimensional saturation unanswered. There are no known examples of $3$-dimensional patterns with $\sat(n,P) = O(1)$ or $\sat(n,P) = \Theta(n)$. It is not even known if, for every multidimensional pattern $P$, there exists an integer $k$ such that $\sat(n,P) = \Theta(n^k)$. 

The argument in Theorem~\ref{thm:2x4s} seems generalizable to most $2 \times n$ matrices. Indeed, it suggests a general strategy of studying patterns of arbitrary size, by analyzing possible witness graphs for them to efficiently search for witnesses.

In Section~\ref{sec:lin_prog}, we presented a linear program for computing the saturation functions of arbitrary 0-1 matrices. We used the linear program to evaluate the saturation functions of some forbidden 0-1 matrices for small values of $m$ and $n$. Based on the results, we made a conjecture on the exact value of the saturation function for some triangular matrices. It would be interesting to compute the exact value of the saturation function for other forbidden 0-1 matrices, using linear programs or other computational methods. Similar methods could also be applied to the problem of computing the exact value of the saturation function of forbidden sequences \cite{anand}.

The witness construction in Proposition~\ref{prop:intermediaries_1} involves starting with a simple pattern and witness, and then inserting additional rows to the pattern and correspondingly multiplying the rows of the witness. This closely resembles the approach used in several of our results around adding empty columns and rows, for example Proposition~\ref{prop:empty_col_h}. It seems likely that a generalization of this construction may be more widely applicable to other patterns, with other constraints on the rows or columns added. In particular, a construction of this kind might be useful in addressing at least parts of Conjectures~\ref{conj:kron} and \ref{conj:empty_rows_cols}, noting that Proposition~\ref{prop:empty_col_h} already covers a special case of the former and is closely related to the latter.

Lastly, we mention the smallest patterns with still unknown asymptotic saturation functions. Several $3 \times 4$ patterns have unclassified $\sat(m_0,n,P)$, and several $5 \times 4$ patterns have unclassified $\sat(n,P)$.

\section*{Acknowledgements}

We thank the organizers of MIT PRIMES-USA, and in particular Dr.~Tanya Khovanova, Dr.~Slava Gerovitch, and Prof.~Pavel Etingof, for making this research possible.


\newpage
\begin{thebibliography}{1}\footnotesize
\bibitem{anand} Anand, J. Geneson, S. Kaustav, and S. Tsai. Sequence saturation. Discrete Appl. Math. 360 (2025) 382-393.
\bibitem{berendsohn21} B.A. Berendsohn. Matrix patterns with bounded saturation function (2020)  \url{https://arxiv.org/abs/2012.14717}. 
\bibitem{berendsohn23} B.A. Berendsohn. An exact characterization of saturation for permutation matrices. Combinatorial Theory 3: 1-35, 2023.
\bibitem{bienstock} D. Bienstock, E. Gy\"ori, An extremal problem on sparse 0-1 matrices, SIAM J. Discrete Math. 4 (1991) 17–27.
\bibitem{BC} R.A. Brualdi and L. Cao. Pattern-avoiding (0, 1)-matrices and bases of permutation matrices. Discrete Appl. Math. 304: 196–211, 2021.
\bibitem{code} A. Duan. Matrix saturation code (2024), \url{https://github.com/TheCoolDinosuar/matrix-saturation}.
\bibitem{Crowdmath2018} P. A. CrowdMath. Bounds on parameters of minimally nonlinear patterns. The Electronic Journal of Combinatorics 25: P1.5, 2018.
\bibitem{EM} P.\ Erd\H os and L.\ Moser, Problem 11, Canadian Math. Bull. 2 (1959), 43.
\bibitem{fulek} R. Fulek, Linear bound on extremal functions of some forbidden patterns in 0–1 matrices. Discrete Math 309: 1736-1739 (2009)
\bibitem{FK} R. Fulek and B. Keszegh. Saturation problems about forbidden 0-1 submatrices. SIAM J. Discrete Math. 35: 1964–1977, 2021.
\bibitem{ngon} Z. F\H{u}redi. The maximum number of unit distances in a convex n-gon. Journal of Combinatorial Theory Series A 55: 316-320, 1990.
\bibitem{FH} Z. F\H{uredi} and  P. Hajnal. Davenport-Schinzel theory of matrices. Discrete Mathematics 103 (3): 233-251, 1992.
\bibitem{geneson2009} J. Geneson. Extremal functions of forbidden double permutation matrices. Journal of Combinatorial Theory Series A 116: 1235-1244, 2009.
\bibitem{geneson19} J. Geneson. Forbidden formations in multidimensional 0-1 matrices. Eur. J. Comb. 78: 147-154, 2019.
\bibitem{geneson21} J. Geneson. A generalization of the K\H{o}v\'{a}ri-S\'{o}s-Tur\'{a}n theorem. Integers, 21: 1–12, 2021.
\bibitem{geneson21a} J. Geneson. Almost all permutation matrices have bounded saturation functions. Electron. J. Combin., 28: 2.16, 2021.
\bibitem{GHLNPW} J. Geneson, A. Holmes, X. Liu, D. Neidinger, Y. Pehova, and I. Wass. Ramsey numbers of ordered graphs under graph operations. arXiv preprint arXiv:1902.00259, 2019.
\bibitem{GS15} J. Geneson, L. Shen, Linear bounds on matrix extremal functions using visibility hypergraphs. Discrete Mathematics 338: 2437-2441, 2015. 
\bibitem{GT17} J. Geneson, P. Tian, Extremal functions of forbidden multidimensional matrices. Discrete Mathematics 340: 2769-2781, 2017.
\bibitem{GTT} J. Geneson, P. Tian, and K. Tung. Formations and generalized Davenport-Schinzel sequences. Integers 22: A111, 2022. 
\bibitem{GT2020} J. Geneson and S. Tsai. Sharper bounds and structural results for minimally nonlinear 0-1 matrices. Electron. J. Combin. 27: 4.24, 2020.
\bibitem{gt23} J. Geneson, S. Tsai, Extremal bounds for pattern avoidance in multidimensional 0-1 matrices. Discrete Mathematics 348: 114303 (2025)
\bibitem{keszegh} B. Keszegh. On linear forbidden submatrices. Journal of Combinatorial Theory, Series A 116: 232-241, 2009.
\bibitem{klazar} M. Klazar, The F\H{u}redi-Hajnal conjecture implies the Stanley-Wilf conjecture. In Formal power series and algebraic combinatorics, pages 250-255. Springer, 2000.
\bibitem{KM} M. Klazar and A. Marcus, Extensions of the linear bound in the Furedi-Hajnal conjecture, Advances in Applied Mathematics 38: 258-266, 2007.
\bibitem{MT} A. Marcus and G. Tardos. Excluded permutation matrices and the Stanley-Wilf conjecture. Journal of Combinatorial Theory Series A 107: 153-160, 2004.
\bibitem{mitchell} J. Mitchell. L1 shortest paths among polygonal obstacles in the plane. Algorithmica 8: 55-88, 1992.
\bibitem{pettie11} S. Pettie, Generalized Davenport-Schinzel sequences and their 0-1 matrix counterparts, J. Combin. Theory Ser. A 118: 1863-1895, 2011.
\bibitem{tardos} G\'{a}bor Tardos. On 0-1 matrices and small excluded submatrices. J. Combin. Theory Ser. A, 111: 266–288, 2005.
\bibitem{tsai23} S. Tsai. Saturation of multidimensional 0-1 matrices. Discrete Math. Lett. 11: 91–95, 2023.
\end{thebibliography}
\end{document}